\documentclass[11pt]{article}

\usepackage{latexsym}
\usepackage{amsmath}
\usepackage{amssymb}
\usepackage{amsfonts}
\usepackage{graphics,graphicx}
\usepackage{color}
\usepackage{comment}
\usepackage{url}
\usepackage{mathrsfs}
\usepackage{latexsym}

\def\claim#1{\begin{trivlist}\item[\hskip\labelsep\bf#1]\it}
\def\endclaim{\end{trivlist}}
\def\ang<#1>{\langle #1 \rangle}

\def\qed{\hfill $\Box$}

\numberwithin{equation}{section}
\newtheorem{theorem}{Theorem}[section]
\newtheorem{lemma}[theorem]{Lemma}

\newtheorem{remark}[theorem]{Remark}
\newtheorem{assumption}{Assumption}[section]
\newcommand{\eproof}{{\mbox{\ }~\hfill
\mbox{\large $\Box$} \par \vskip 10pt}}
\newcommand{\pf}{\noindent{\bf Proof}}
\newcommand{\Om}{\Omega}

\date{}

\begin{document}


\title{Uniform Decaying Property of Solutions for Anisotropic Viscoelastic Systems}

\author{Maarten V. de Hoop \thanks{Departments of Computational and Applied
Mathematics and Earth, Environmental and Planetary Sciences, Rice University, Houston TX, USA, (Email: mdehoop@rice.edu)}
\qquad Ching-Lung Lin\thanks{Department of Mathematics, National Cheng-
Kung University, Tainan 701, Taiwan. Partially supported by the
Ministry of Science and Technology of Taiwan. (Email:
cllin2@mail.ncku.edu.tw)}\qquad Gen Nakamura\thanks{Department of
Mathematics, Hokkaido University, Sapporo 060-0808, Japan and Research Center of Mathematics for Social Creativity, Research Institute for Electronic Science, Hokkaido University, Sapporo 060-0812, Japan. Partially supported by JSPS KAKENHI (Grant No. JP19K03554 and 22K03366). (Email: gnaka@math.sci.hokudai.ac.jp)
}}

\maketitle

\begin{abstract}
The paper concerns about the uniform decaying property (abbreviated by UDP) of solutions for an anisotropic viscoelastic system in the form of integrodifferential system (abbreviated by VID system) with mixed type boundary condition. The mixed type condition consists of the homogeneous displacement boundary condition and a homogeneous traction boundary condition or with a dissipation. By using a dissipative structure of this system, we will prove the UDP in a unified way for the two cases, which are, when the time derivative of relaxation tensor decays with polynomial order and it decays with exponential order, by clarifying and modifying the argument of Rivera-Lapa \cite{R-L}.
\end{abstract}

\renewcommand{\theequation}{\thesection.\arabic{equation}}

\section{Introduction}\label{sec1}

In this paper, we study the uniform decaying property (abbreviated by UDP) of solutions to the initial boundary value problem with mixed type condition for a general anisotropic viscoelastic system containing a memory term (abbreviated by VID system) by clarifying and modifying the argument of Rivera-Lapa \cite{R-L} in a unified and systematic way, so that the argument and its idea can be understood more easily. The mixed type boundary condition consists of the homogeneous displacement condition and the traction free boundary condition. We refer this initial boundary value problem as (MP). We will prove the UDP for the two cases, which are, when the time derivative of relaxation tensor decays with polynomial order in time and it decays with exponential order in time. Corresponding to these cases, solutions of the VID system decay with polynomial order in time and exponential order in time, respectively.

In order to formulate the initial boundary value problem, let $\Omega\subset\mathbb{R}^{n}$ with $n\in\mathbb{N}$ be a bounded domain. Denote its boundary by $\partial\Omega$. For $n\ge2$, we assume that $\partial\Omega$ is $C^1$ smooth and connected, and for $n=1$, we assume that it consists of two points.
Divide $\partial\Omega$ into $\partial\Omega = \overline{\Gamma_D}\cup\overline{\Gamma_N}$, where $\Gamma_D,\,\Gamma_N\subset\partial\Omega$ are open and assume that $\Gamma_D\not=\emptyset$, $\Gamma_D\cap\Gamma_N=\emptyset$. If $n\ge 3$, then
we further assume that their boundaries $\partial\Gamma_D,\,\partial\Gamma_N$ are  Lipschitz smooth.

We consider the following initial boundary value problem (MP)'
\begin{equation}\label{mixed prob}
\left\{
\begin{array}{rcll}
\rho\partial_t^2 u(\cdot, t)&=&\nabla\cdot\{C(\cdot)\nabla u(\cdot,t)-\displaystyle{\int_0^t} &\!\!\!G(\cdot,t-\tau)\nabla u(\cdot,\tau)\,d\tau\}\qquad(t>0),
\\[0.25cm]
u&=&0\,\,\mbox{\rm on}\,\,\Gamma_D\times(0,\infty),\,\,&\mathcal{T}u+s\partial_t u=0\,\,\mbox{\rm on}\,\,\Gamma_N\times(0,\infty),
\\[0.25cm]
u&=&f_1\in H_+,\,\,&\partial_t u=f_2\in  L^2(\Omega)\,\,\mbox{\rm on}\,\,\Omega\times\{0\},
\end{array}
\right.
\end{equation}
where $s\geq 0$ is a non-negative constant and $\mathcal{T}u$ is the traction given by
$$
\mathcal{T}u(\cdot,t)=\big(C(\cdot)\nabla u(\cdot,t)-\int_0^t\,G(\cdot,t-\tau)\nabla u(\cdot,\tau)\,d\tau\big)\nu
$$
with the unit outer normal vector $\nu$ of $\partial\Omega$ and $H_+:=\{f_1\in H^1(\Omega): f_1\big|_{\Gamma_D}=0\}$.
Here $0<\rho_0\le\rho\in L^\infty(\Omega)$ with a positive constant $\rho_0$, $C=(C_{ijk\ell})$ and $G=(G_{ijk\ell})$ denote the elasticity tensor and $-1$ times the time derivative of the relaxation tensor, respectively. The VID system is the first equation of \eqref{mixed prob} and the aforementioned (MP) is given by putting $s=0$ for (MP)'. If we move the dissipation term $\partial_t u$ to the right hand side and consider $-s\partial_t u$ as a control, we have a boundary controllability problem.

There are several ways to derive the VID system, which is a system of equations in \eqref{mixed prob}. Its abstract derivation is given in \cite{C}. The well known Maxwell model and its extended version yield the VID system via some system describing the deformations of springs and dashpots. Here, the extension includes the higher space dimensional case and the case that the media for the model are anisotropic. We will call such a system a $\text{\rm spring-dashpot}$ system. The standard linear solid (SLS) model and its extended version are  a little bit more general cases of the extended Maxwell model; see \cite{DKLN} for further details on their full tensorial realizations. Recently, such realization has been developed for the extended Burgers model (see \cite{EBM}). In space-dimension $1$, when the system reduces to a scalar equation and anisotropy no longer comes into play, the VID form of the Burgers model is given in Appendix A. Concerning the method of deriving the VID system from the extended spring-dashpot system, we note that an approach proposed by Carcione \cite{Carcione} using spectral representation for analyzing wave fields of the VID system is also useful (see \cite{EBM}).

Before moving to a new subsection, it should be noted that in this paper, we are primarily interested in the anisotropic case for $C,\,G$, and will never specialize our arguments to the isotropic case.

\subsection{Main results}

We first state assumptions (Assumptions \ref{polynomial}) for the polynomial order decay of solutions and assumptions (Assumptions \ref{exponential}) for the exponential order decay of solutions of (MP) and (MP)'.
Henceforth, we assume $\rho=1$ for simplicity of presentation. Then, the assumptions are as follow.
\begin{assumption}(polynomial order decay)\label{polynomial}${}$
\newline
(i)
$C\in L^\infty(\Omega)$ and $G(\cdot,t)\in C^2([0,\infty); L^\infty(\Omega))$.
There exist constants  $\kappa_1>0$, $\kappa_2>0$, $\kappa_3>0$ and $p>2$ such that for any $\text{a.e.}\,x\in\Omega,\,\,t\in[0,\infty)$,
\begin{equation}\label{1.2}
 -\kappa_1G(x,t)\le\dot{G}(x,t)\le-\kappa_2G^{1+\frac{1}{p}}(x,t),\,\, \ddot G(x,t)\le\kappa_3G^{1+\frac{1}{p}}(x,t),
\end{equation}
where $\dot G(x,t)=\partial_t G(x,t)$ and $\ddot G(x,t)=\partial_t^2 G(x,t)$, respectively. Also, there exists a constant $\kappa_4>0$ such that
\begin{equation}\label{1.3}
|\dot{G}(x,t)|+|G(x,t)|\leq\kappa_4(1+t)^{-p},\,\,\text{a.e $x\in\Omega$}, \,\,t\ge0.
\end{equation}

\medskip
\noindent
(ii) (major symmetry) $C_{ijk\ell}=C_{k\ell ij},\,\,G_{ijk\ell}=G_{k\ell ij}\,\,\text{a.e. in $\Omega$}\,\,,t\ge0,\,\,1\le i,j,k,\ell\le n$.

\medskip
\noindent
(iii) (strong convexity) There
exist constants $\alpha_0>0$ and $\beta_0>0$ such that for any $n\times n$ symmetric matrix $w=(w_{ij})$
\begin{equation}\label{1.4}
\alpha_0|w|^2\le(Cw):w\le\beta_0|w|^2,\,\,
\alpha_0|w|^2\le(Gw):w\le\beta_0|w|^2\,\,\text{a.e. in $\Omega$},\,\,t\ge0
\end{equation}
where the notation ``$:$'' is defined as $(Cw):w=\sum_{i,j,k,\ell=1}^n C_{ijk\ell}w_{ij}w_{k\ell}$. The lower estimates in \eqref{1.4} are referred to as the strong convexity condition for $C,\,G$.

\medskip
\noindent
(iv) There exist constants $\mu_0>0,\,\nu_0>0$ such that for $n\times n$ symmetric matrix $w=(w_{ij})$,
\begin{equation}\label{1.5}
\mu_0|w|^2\leq\{(C(\cdot)-\int_0^\infty G(\cdot,t)dt)\,w\}:w\le\nu_0|w|^2\,\,\text{a.e. in $\Omega$}. 
\end{equation}
That is $C(\cdot)-\int_0^\infty G(\cdot,t)dt$ satisfies the strong convexity condition.
\end{assumption}

\begin{assumption}(exponential order decay)\label{exponential} 
\begin{itemize}
\item[{\rm (i)'}] $C\in L^\infty(\Omega)$ and $G(\cdot,t)\in C^2([0,\infty); L^\infty(\Omega))$.
There exist constants  $\kappa_j>0,\,j=1,\cdots,4$ and $\tilde\kappa_4>0$ such that for any $\text{a.e.}\,x\in\Omega,\,\,t\in[0,\infty)$,
\begin{equation}\label{derivatives}
-\kappa_1 G(x,t)\le\dot G(x,t)\le-\kappa_2 G(x,t),\,\ddot G(x,t)\le\kappa_3 G(x,t)
\end{equation}
and
\begin{equation}\label{exp decay of the |G| and |dot G|}
|G(x,t)|+|\dot G(x,t)|\le \kappa_4e^{-\tilde{\kappa}_4 t}.
\end{equation}
\item[{\rm (ii)'}] $C(x),\,G(x,t)$ satisfy the major symmetry condition.
\item[{\rm (iii)'}]
$C(x),\,G(x,t)$ satisfy
the strong convexity condition.
\item[{\rm (iv)'}] $C(x)-\int_0^\infty G(x,t)\,dt$ satisfies the strong convexity condition.
\end{itemize}
\end{assumption}

Our main results are as follow.

\begin{theorem}\label{mainresult}
For any given initial data $(f_1,f_2)\in H_+\times H_+$, the solution
$u\in C^2([0,\infty); H^1(\Omega))$ of \eqref{mixed prob} with the property $\partial_t^3 u\in C([0,T]; L^2(\Omega))$, $\partial_t^4 u\in C([0,\infty); H_-)$ converges to zero at polynomial rate in time $t$ as $t\rightarrow\infty$ and this rate does not depend on $u$. Here, $H_-$ denotes the dual space of $H_+$.
\end{theorem}

\begin{theorem}\label{mainresult2}
For any given initial data $(f_1,f_2)\in H_+\times   H_+$, the solution
$u\in C^2([0,\infty); H^1(\Omega))$ of \eqref{mixed prob} with the property $\partial_t^3 u\in C([0,T]; L^2(\Omega))$, $\partial_t^4 u\in C([0,\infty); H_-)$ converges to zero at exponential rate in time $t$ as $t\rightarrow\infty$ and this rate does not depend on $u$.
\end{theorem}

\begin{remark}\label{several remarks}${}$
\begin{itemize}
\item[(i)] The UDP for (MP) with the polynomial decay rate and the exponential decay rate follow by putting $s=0$ for Theorem \ref{mainresult} and Theorem \ref{mainresult2}, respectively.
\item[(ii)](unique existence of solutions) Suppose the initial data satisfies the so called the compatibility condition of order 2 (i.e. $(u_0,u_1,0)\in B_2(T)$ in \cite{D1} with the initial data $(u_0,u_1)$ and the inhomogeneous term $0$ of the VID system). Then, for any $T>0$, the unique existence of the solution $u\in C^2([0,T]; H^1(\Omega))$ of \eqref{mixed prob} over $[0,T]$ with $\partial_t^3 u\in C([0,T]; L^2(\Omega))$ follows by mimicking the proof of Theorem~2.2 in \cite{D1} with a modification and the proof of Theorem~2.3 therein. The modification is summarized in \ref{App:modification}.
Also, since $T > 0$ is arbitrary, we immediately obtain the unique existence of a solution $u \in C^2([0,\infty); H^1(\Omega))$ of \eqref{mixed prob}.
\item[(iii)] (lower bounds) It is important to assume that there are lower bounds for $\dot G$ given in \eqref{1.2} and \eqref{derivatives} for the polynomially decaying kernel $G(x,t)$ and exponentially decaying kernel, respectively. For the polynomially decaying kernel, this lower bound was not stated in \cite{R-L}.
\item[(iv)] (separable kernel)
Let each positive full symmetric $4$ rank tensor $\mathscr{G}_j(x),\,\hat G_j(x)\in L^\infty(\Omega),\,1\le j\le J$ satisfies
$$
\mathscr{G}_j(x)\ge g I,\,\hat G_j(x)\ge g I\,\,\text{a.e.}\,\,x\in\Omega,
$$
where $g>0$ is a constant and $I$ is the $4$ rank identity tensor. As a generalization of seperable kernel $e^{-\mathscr{G}t}\hat G(x)$ with positive full symmetric $4$ rank tensors $\mathscr{G}$ and $\hat{G}(x)\in L^\infty(\Omega)$ satisfying
$$
\mathscr{G}\ge g I,\,\hat G(x)\ge g I\,\,\text{a.e.}\,\,x\in\Omega,
$$
we consider
$$
G(x,t):=\displaystyle\sum_{j=1}^J e^{-\mathscr{G}_j(x)\, t}\hat G_j(x)
$$
and call it separable. Here, the above $4$ rank full symmetric tensors are considered as either tensor block symmetric matrices (see \cite{EBM}) or in terms of the Voigt notation. Also, $e^{-\mathscr{G}_j(x)t}$ is understood as $\sum_{n=0}^\infty (n!)^{-1} t^n(-1)^n\mathscr{G}_j(x)^n$. It is easy to see that (iv)' of Assumption \ref{exponential} becomes
$$
(\ast)\, C(x)-\displaystyle\sum_{j=1}^J\mathscr{G}_j(x)^{-1}\hat{G}_j(x)\,\text{is strongly convex}.
$$
All other assumptions of Assumption \ref{exponential} are satisfied by the assumptions we made for the separable kernel $G(x,t)$.

\item[(v)] (examples) Well known spring-dashpot models in mechanics such as the Maxwell model, standard linear solid (abbreviated by SLS) model, Burgers model and their higher space dimensional and anisotropic extended version yield VID systems under some condition at $t=0$. $G(x,t)$ of these VID systems can be separable. For more details, see  \ref{secA_M} and \cite{EBM}.
\end{itemize}
\end{remark}

\subsection{Prior works}

There are several studies on the asymptotic behavior of solutions of VID system as follows. Some abstract schemes for an
integrodifferential equation were developed for the well-posedness of the abstract Cauchy problem and the asymptotic behavior of solutions in \cite{D1} and \cite{D2}, respectively. It is shown in \cite{D2} that solutions of VID system satisfying the Dirichlet boundary condition decay to zero as the time tends to infinity. Concerning the decay rate of the solutions, there are two kinds of decay rates. One is the polynomial order decay rate and the other is the exponential order decay rate. As for a polynomial order decay rate, an energy method given in \cite{R-L} was really a big break through. In this paper, polynomial order decay of solutions satisfying the Dirichlet boundary condition was given for very general anisotropic viscoelastic systems by using an energy norm given in \cite{D1}.
To our knowledge, the earliest source for this norm is in \cite{D1}. Ahead of discussing more about its argument, we give an important results for the exponential order decay rate. By looking at the structure of the resolvent of the elliptic system obtained by taking the Laplace transform of VID system with respect to $t$, the first result on the exponential decay of solutions of the initial boundary value problem for VID system with the Dirichlet boundary condition was given in \cite{Fabrizio} under some assumptions. Also, for a special isotropic viscoelastic integrodifferential system with exponentially decaying relaxation tensor, the exponential decay of solutions satisfying the Dirichlet boundary condition was given in \cite{R} by using the aforementioned energy norm.
Then in \cite{NO}, it was further generalized to the VID system with a mixed type boundary condition in a way along the lines of the argument in \cite{R, R-L}.

Here, we would like to discuss more about the argument of \cite{R-L}. In the paper, if the initial displacement is not zero, we could not figure out how to handle the last term of the second identity in Lemma 2.3 of \cite{R-L} to have the estimate in Remark 2.1. In our paper, the third term of (3.3) of Lemma 3.2 is exactly the same as the above mentioned last term. Since this term contains the initial information, we amended the energy type function $\mathcal{L}(t,u)$ in \cite{R-L} to contain the suitable initial information in \eqref{mathcal L}. This amendment enabled us to obtain the estimate \eqref{3.17} in Lemma \ref{ident5} which is similar to Remark 2.1 of \cite{R-L}. Then, following the argument in \cite{R-L} with necessary modifications, we showed the UDP with polynomial order decay rate for solutions of \eqref{mixed prob}. By this clarification and modification, we could give two kinds of the UDP with different decay rates by a unified argument.

\subsection{Novelty of the results}
The novelty of the results in this paper are as follow.
\begin{itemize}
\item [(i)] We extended the polynomial order decay result of \cite{R} to the initial boundary value problem for the VID system with mixed type boundary condition consisting of homogeneous displacement condition and the traction free boundary condition. We further extended the result to the case the mentioned traction free boundary condition, i.e. $\mathcal{T}u=0$, is replaced by  $\mathcal{T}u+s\partial_t u=0$.
\item[(ii)] We extended the exponential order decay result of \cite{NO} which was given under zero initial displacement and $G(x,t)=e^{-\kappa t}\hat G(x)$ with a constant $\kappa>0$ to the general initial data and general $G(x,t)$ given in Assumption \ref{exponential}. Likewise (i), we further extended the result to the case the mentioned traction free boundary condition, i.e. $\mathcal{T}u=0$, is replaced by  $\mathcal{T}u+s\partial_t u=0$.
\end{itemize}

\subsection{Organization of this paper}
The rest of this paper is organized as follows. The main part of the paper is devoted to the proof of Theorem~\ref{mainresult}. The idea of proof relies on that of \cite{R-L}. By clarifying and refining the idea given in \cite{R-L}, we provide the proof of Theorem~\ref{mainresult} which is very systematic so that the proof of Theorem~\ref{mainresult2} can be given almost the same as that of Theorem~\ref{mainresult}. The outline of each section is given as follows. In Section 2, we introduce some notations and give the strategy of proof. Then we provide some basic identities and inequalities in Section 3. Having these, we carefully carry out the strategy in Section 4 to prove Theorem~\ref{mainresult}. In Section 5, we give a proof of Theorem \ref{mainresult2}. As we mentioned before, the proof of Theorem~\ref{mainresult2} is almost the same as that of Theorem \ref{mainresult}, we only give necessary modifications of the proof of Theorem \ref{mainresult}. Appendix A gives several examples of VID associated to extended spring-dashpot model that satisfy the conditions of Assumption \ref{polynomial} and Assumption \ref{exponential}. Concerning the unique solvability of VID system, Appendix B provides modifications of the proof of Theorem~2.2 in \cite{D1} that is necessary for our mixed type boundary condition.

\section{Outline of the proof of Theorem~\ref{mainresult} and preliminaries}

To begin with, unless otherwise specified, we have to remark here that $u\in C^2([0,\infty); H^1(\Omega))$ in this paper is the solution to \eqref{mixed prob} with the property $\partial_t^3 u\in C([0,T]; L^2(\Omega))$, $\partial_t^4 u\in C([0,\infty); H_-)$. Knowing this, let
\[
F(\cdot,t) := \gamma G(\cdot,t)+\dot{G}(\cdot,t)
\]
and
introduce the quantities
\begin{multline}
E(t,u):=\frac12\big[\int_\Omega|\dot{u}(\cdot,t)|^2+(G\Box\partial u)(\cdot,t)\,dx
\\
+\int_\Omega\{(C(\cdot)-\int_0^tG(\cdot,\tau)\,d\tau)\nabla u(\cdot,t)\}:\nabla u(\cdot,t)\,dx\big],
\end{multline}
\begin{equation}
G\Box\partial u(\cdot,t):=\int_0^t\big\{G(\cdot,t-\tau)\nabla\big(u(\cdot,t)-u(\cdot,\tau)\big)\big\}:\nabla\big(u(\cdot,t)-u(\cdot,\tau)\big)\,d\tau,
\end{equation}
\begin{multline}
K(t,u):=\frac12\int_\Omega|\ddot{u}|^2\,dx+\frac12\int_\Omega\big(C(\cdot)\nabla\dot{u}(\cdot,t)\big):\nabla\dot{u}(\cdot,t)\,dx
\\
-\int_\Omega\big(G(\cdot,0)\nabla u(\cdot,t)\big):\nabla\dot{u}(\cdot,t)\,dx+
\gamma\int_\Omega\big(C(\cdot)\nabla u(\cdot,t)\big):\nabla\dot{u}(\cdot,t)\,dx
\\
-\int_\Omega\big(\int_0^t F(\cdot,t-\tau)\nabla u(\cdot,\tau)\,d\tau\big):\nabla\dot{u}(\cdot,t)\,dx,
\end{multline}
\begin{multline}
I(t,u):=\int_\Omega\ddot{u}(\cdot,t)\dot{u}(\cdot,t)\,dx-\frac12\int_\Omega\big(G(\cdot,0)\nabla u(\cdot,t)\big):\nabla u(\cdot,t)\,dx,
\\
-\frac12\int_\Omega\big(\int_0^t\dot{G}(\cdot,\tau)\,d\tau\nabla u(\cdot,t)\big):\nabla u(\cdot,t)\,dx+\frac12\int_\Omega\dot{G}\Box\partial u(\cdot,t)\,dx,
\end{multline}
Now, we define
\begin{equation*}
\begin{aligned}
B(t,u):= N_3E(t,u)+E(t,\dot u)-\int_\Omega\big(G(\cdot,t)\nabla u(\cdot,0)):\nabla\dot{u}(\cdot,t)\,dx,
\end{aligned}
\end{equation*}
where $N_3$ is large enough such that
\begin{equation}\label{B}
\begin{cases}
\begin{aligned}
&\frac{1}{2}N_3 E(0,u)+\frac{1}{2}E(0,\dot u)\leq B(0,u)\leq 2N_3E(0,u)+2E(0,\dot u)\\
&\frac{1}{2}N_3 E(t,u)+\frac{1}{2}E(t,\dot u)\leq B(t,u)+\hat{c}_\delta(1+t)^{-p}E(0,u).
\end{aligned}
\end{cases}
\end{equation}
In order to see that we can take $N_3$ to have \eqref{B}, we note that for any $\delta>0$, there exist $c_\delta>0,\,\hat c_\delta>0$ such that
\begin{equation*}
\begin{aligned}
&|\int_\Omega\big(G(\cdot,t)\nabla u(\cdot,0)):\nabla\dot{u}(\cdot,t)\,dx|\\
\leq & \delta \int_\Omega|\nabla\dot{u}(\cdot,t)|^2\,dx+c_\delta (1+t)^{-p}\int_\Omega|\nabla u(\cdot,0)|^2\,dx\\
\leq & \frac{1}{2} E(t,\dot u)+\hat{c}_\delta (1+t)^{-p}E(0,u).
\end{aligned}
\end{equation*}
We also let
\begin{equation}\label{mathcal L}
\mathcal{L}(t,u):=N_1 [ B(t,u)+2\hat{c}_\delta(1+t)^{-p}E(0,u)]+K(t,u)+(\gamma-\theta)I(t,u)+\omega\int_\Omega\dot{u}u\,dx,
\end{equation}
where $N_1,\,\gamma,\,\theta,\,\omega$ are positive constants which will satisfy some condition given later. Here and hereafter, whenever we have $wz$ for any two vectors $w,\,z$, it means the real inner product of $w$ and $z$.

We remark here that the first two terms of the right hand side of \eqref{mathcal L} are the most important quantities and the other terms are necessary to estimate when we differentiate the first two terms with respect to $t$.

We will estimate $\mathcal{L}(t,u)$ from below by a positive constant times the sum
$\int_\Omega|\nabla u(\cdot,t)|^2\,dx$ with some other positive terms depending on $u$ and $\frac{d}{dt} \mathcal{L}(t,u)$ from above. More precisely, we will prove in Lemma \ref{lem4.2} and Lemma \ref{lem4.3} that
\begin{equation}\label{2.1}
\mathcal{L}(t,u)\ge M_2\int_\Omega\big\{|\ddot{u}|^2+|\dot{u}|^2+|\nabla\dot{u}|^2+|\nabla u|^2\big\}\,dx,\,\,t>0
\end{equation}
and
\begin{equation}\label{2.2}
\frac{d}{dt}\mathcal{L}(t,u)\leq M_3(1+t)^{-p_m-1}-M_2\mathcal{L}(t,u)^{(1+\frac{1}{p_m})},\,\,t>0
\end{equation}
for some positive constants $M_1,\,M_2,\,M_3$ and $m$ is a positive integer such that $p_m=2^m-1<p-1$.

\medskip
To estimate $\mathcal{L}(t,u)$ from above, we prepare the following lemma.

\begin{lemma}\label{lem2.1}
Let $M_2,\,M_3>0$ and $q>2$.  Assume that $y\in C^1[0,\infty)$ and $y(t)\geq 0,\,t\geq 0$ satisfying
\begin{equation}\label{2.3}
\begin{aligned}
\frac{d}{dt}y(t)\leq-M_2y(t)^{(1+\frac{1}{q})}+M_3(1+t)^{-q-1},\,\,t>0,
\end{aligned}
\end{equation}
then we have that
\begin{equation}\label{2.4}
\begin{aligned}
y(t)\leq q^q\big[(y(0)+2M_2^qM_3)^{-\frac{1}{q}}+(2M_2)^{-1}M_3^{-\frac{1}{q}}t\big]^{-q},\,\,t>0.
\end{aligned}
\end{equation}
\end{lemma}

\pf:\
By denoting $\tilde{y}(t)=M_2^q\,y(t)$, we can assume that the constants  $M_2$ and $M_3$ in \eqref{2.3} become $1$ and $M_4=M_2^qM_3$. We define a new function $z(t):= y(t)+ l(t)$, where $l(t)=2M_4(1+t)^{-q}$. Direct computations give that
\begin{equation}\label{2.5}
\begin{aligned}
\frac{d}{dt}z(t)=&\frac{d}{dt}y(t)+\frac{d}{dt}l(t)\\
\leq & -y(t)^{(1+\frac{1}{q})}+M_4(1+t)^{-q-1}-2qM_4(1+t)^{-q-1}\\
\leq & -y(t)^{(1+\frac{1}{q})}-2M_4(1+t)^{-q-1}\\
=& -y(t)^{(1+\frac{1}{q})}-(2M_4)^{-\frac{1}{q}}\,l(t)^{(1+\frac{1}{q})}\\
\leq & - (2)^{-(1+\frac{2}{q})}M_4^{-\frac{1}{q}}\,z(t)^{(1+\frac{1}{q})}\\
\leq & - M_5\,z(t)^{(1+\frac{1}{q})},
\end{aligned}
\end{equation}
where $M_5=(2M_2)^{-1}M_3^{-\frac{1}{q}}$.
From \eqref{2.5}, we have
\begin{equation}\label{2.6}
\begin{aligned}
\frac{d}{d\tau}\big(-qz(\tau)^{-\frac{1}{q}}\big)\leq  - M_5.
\end{aligned}
\end{equation}
Integrating  this over $[0,t]$, we have
\begin{equation}\label{2.7}
\begin{aligned}
y(t)\leq z(t)\leq  q^q\big[(y(0)+2M_2^q\,M_3)^{-\frac{1}{q}}+(2M_2)^{-1}M_3^{-\frac{1}{q}}t\big]^{-q}.
\end{aligned}
\end{equation}
\qed

\medskip
Now applying this lemma to $\mathcal{L}(t,u_\epsilon)$ with $q=p_m$ and $ u_\varepsilon:=(\mathcal{L}(0,u)+\epsilon)^{-1/2}u$, $\epsilon>0$,
we have
$$
\mathcal{L}(t, u_\epsilon)\le q^q[(\mathcal{L}(0, u_\epsilon)+2M_2^q M_3)^{-\frac1{q}}+(2M_2)^{-1} M_3^{-\frac1{q}} t]^{-q},\,\,t>0.
$$
Then, by $\mathcal{L}(t, u_\epsilon)=(\mathcal{L}(0,u)+\epsilon)^{-1}\mathcal{L}(t,u)$ and

$$\displaystyle\lim_{\epsilon\rightarrow 0}\mathcal{L}(0,u_\epsilon)=
\left\{
\begin{array}{ll}
1\,\,\,\text{if}\,\,\mathcal{L}(0,u)>0,\\
0\,\,\,\text{if}\,\,\mathcal{L}(0,u)=0,
\end{array}
\right.
$$
we have
\begin{equation}\label{decay mathcal L}
\mathcal{L}(t,u)\le \mathcal{L}(0,u)\,a\,(b_1+b_2t)^{-q}
\end{equation}
with $a=q^q$, $b_1=(1+2M_2^q\, M_3)^{-1/q}$ and $b_2=(2M_2)^{-1} M_3^{-1/q}$.

Then to finish the proof of Theorem \ref{mainresult}, we only need to note \eqref{B}, \eqref{2.1} and Poincar\'e's lemma.

\section{Basic identities and inequalities for (MP)'}

To prove estimate (\ref{2.2}) we need a collection of basic identities and inequalities which we are going to give in this section. In deriving these identities and inequalities, we will pay attention to the mixed type boundary condition and the constants of the inequalities in Assumption \ref{polynomial}.

\subsection{Basic identities}

\medskip
\noindent
We first simply cite the following identity given as Lemma 2.1. in \cite{R-L}.
\begin{lemma}\label{ident1}
\begin{equation}\label{3.1}
\begin{aligned}
&\int_\Omega\big(\int_0^t G(\cdot,t-\tau)\nabla v(\cdot,\tau)\,d\tau\big):\nabla\dot{v}(\cdot,t)\,dx\\
=&-\frac12\int_\Omega\big(\frac{d}{dt}G\Box\partial v\big)(\cdot,t)\,dx+\frac12\int_\Omega\big(\dot{G}\Box\partial v\big)(\cdot,t)\,dx\\
&+\frac12\int_\Omega\frac{d}{dt}\big(\int_0^tG(\cdot,\tau)\,d\tau\,\nabla v(\cdot,t)\big):\nabla v(\cdot,t)\,dx\\
&-\frac{1}{2}\int_\Omega\big(G(\cdot,t)\nabla v(\cdot,t)\big):\nabla v(\cdot,t)\,dx,\,\,\,t>0.
\end{aligned}
\end{equation}
\end{lemma}

\begin{lemma}\label{ident2}
\begin{equation}\label{3.2}
\begin{aligned}
\frac{d}{dt} E(t,u)=&-\frac{1}{2}\int_\Omega\big(G(\cdot,t)\nabla u(\cdot,t)\big):\nabla u(\cdot,t)\,dx+\frac12\int_\Omega\dot{G}\Box\partial u(\cdot,t)\,dx\\
&-s\int_{\Gamma_N}|\dot{u}(\cdot,t)|^2\,d\sigma,
\end{aligned}
\end{equation}
\begin{equation}\label{3.3}
\begin{aligned}
\frac{d}{dt} E(t,\dot{u})=&-\frac{1}{2}\int_\Omega\big(G(\cdot,t)\nabla \dot{u}(\cdot,t)\big):\nabla \dot{u}(\cdot,t)\,dx+\frac12\int_\Omega\dot{G}\Box\partial \dot{u}(\cdot,t)\,dx\\
&+\int_\Omega\big(G(\cdot,t)\nabla u(\cdot,0)):\nabla\ddot{u}(\cdot,t)\,dx-s\int_{\Gamma_N}|\ddot{u}(\cdot,t)|^2\,d\sigma,\,\,t>0.
\end{aligned}
\end{equation}
\end{lemma}

\pf.
Let us multiply the viscoelastic equation in (\ref{mixed prob}) by $\dot{u}(\cdot,x)$ to have
\begin{equation}\label{3.4}
\begin{aligned}
&\frac{1}{2}\frac{d}{dt}\Big\{\int_\Omega|\dot{u}(\cdot,t)|^2\,+\big(C(\cdot)\nabla u(\cdot,t)\big):\nabla u(\cdot,t)\,dx \Big\}\\
=&\int_\Omega\big(\int_0^t G(\cdot,t-\tau)\nabla u(\cdot,\tau)\,d\tau\big):\nabla\dot{u}(\cdot,t)\,dx-s\int_{\Gamma_N}|\dot{u}(\cdot,t)|^2\,d\sigma.
\end{aligned}
\end{equation}
Using Lemma \ref{ident1}, our first assertion holds. To show the second identity, we take the time derivative of the viscoelastic equation in (\ref{mixed prob}) so that
\begin{equation}\label{3.5}
\begin{aligned}
\stackrel{(3)}{u}(\cdot,t)+\nabla\cdot\Big\{-C(\cdot)\nabla\dot{u}(\cdot,t)\,+G(\cdot,0)\nabla u(\cdot,t)\,+\int\limits_{0}^{t}\dot{G}(\cdot,t-\tau)\nabla u(\cdot,\tau)\,d\tau \Big\}=0,
\end{aligned}
\end{equation}
where $\stackrel{(3)}{u}(\cdot,t)$ denotes the third order derivtive on $u(\cdot,t)$ with respect to $t$.
Multiplying this by $\ddot{u}(\cdot,t)$ and using again Lemma \ref{ident1}, we have the second identity
\begin{equation}\label{3.6}
\begin{aligned}
\frac{1}{2}\frac{d}{dt}\Big\{\int_\Omega|\ddot{u}(\cdot,t)|^2\,+\big(C(\cdot)\nabla \dot{u}(\cdot,t)\big):\nabla \dot{u}(\cdot,t)\,dx \Big\}=I_1+I_2,
\end{aligned}
\end{equation}
where
$$ I_1=\int_\Omega\big(\int_0^t G(\cdot,t-\tau)\nabla u(\cdot,\tau)\,d\tau\big):\nabla\dot{u}(\cdot,t)\,dx,$$
$$I_2=\int_\Omega\big(G(\cdot,t)\nabla u(\cdot,0)):\nabla\ddot{u}(\cdot,t)\,dx-s\int_{\Gamma_N}|\ddot{u}(\cdot,t)|^2\,d\sigma.$$
\eproof
\begin{lemma}\label{ident3}
\begin{equation}\label{3.7}
\begin{aligned}
&\frac{d}{dt}\big\{K(t,u)+(\gamma-\theta)I(t,u)\big\}\\
=&-\theta\int_\Omega|\ddot{u}(\cdot,t)|^2\,dx+\theta\int_\Omega\big(C(\cdot)\nabla\dot{u}(\cdot)\big):\nabla\dot{u}(\cdot,t)\,dx\\
&-\int_\Omega\big(G(\cdot,0)\nabla\dot{u}(\cdot,t)\big):\nabla\dot{u}(\cdot,t)\,dx\\
&-\frac{1}{2}(\gamma-\theta)\int_\Omega\big(\dot{G}(\cdot,t)\nabla u(\cdot,t)\big):\nabla u(\cdot,t)\,dx\\
&+\frac12(\gamma-\theta)\int_\Omega\ddot{G}\Box\partial u(\cdot,t)\,dx
-\int_\Omega\big(F(\cdot,t)\nabla u(\cdot,t)\big):\nabla\dot{u}(\cdot,t)\,dx\\
&+\int_\Omega\big\{\int_0^t\dot{F}(\cdot,t-\tau)\nabla\big(u(\cdot,t)-u(\cdot,\tau)\big)\,d\tau\big\}:\nabla\dot{u}(\cdot,t)\,dx\\
&-s\int_{\Gamma_N}|\ddot{u}(\cdot,t)|^2\,d\sigma+s(\theta-2\gamma)\int_{\Gamma_N}\dot{u}(\cdot,t)\ddot{u}(\cdot,t)\,d\sigma,\,\,\,t>0.
\end{aligned}
\end{equation}
\end{lemma}
\pf.
First, we sum the viscoelastic equation multiplied by $\gamma$ and the time derivative of the viscoelastic equation in (\ref{mixed prob}) to obtain
\begin{equation}\label{3.8}
\begin{aligned}
&\stackrel{(3)}{u}(\cdot,t)\,+\gamma \ddot{u}(\cdot,t)+\nabla\cdot\Big\{-C(\cdot)\nabla\dot{u}(\cdot,t)\,+G(\cdot,0)\nabla u(\cdot,t)\Big\}\\
=&-\nabla\cdot\Big\{\int\limits_{0}^{t}F(\cdot,t-\tau)\nabla u(\cdot,\tau)\,d\tau\,+\gamma C(\cdot) \nabla u(\cdot,t) \Big\}
\end{aligned}
\end{equation}
Hence, multiplying by $\ddot{u}(\cdot,t)$ and integrating in $\Omega$, we have
\begin{equation}\label{3.9}
\begin{aligned}
&\frac{1}{2}\frac{d}{dt}\Big\{\int_\Omega|\ddot{u}(\cdot,t)|^2+(C(\cdot)\nabla \dot{u}(\cdot,t)):\nabla \dot{u}(\cdot,t)\,dx\Big\}=J_1+J_2+J_3+J_4,
\end{aligned}
\end{equation}
where
$$J_1=-\gamma\int_\Omega |\ddot{u}(\cdot,t)|^2\,dx+\int_\Omega (G(\cdot,0)\nabla u(\cdot,t)):\nabla\ddot{u}(\cdot,t)\,dx,$$
$$J_2=-\gamma\int_\Omega (C(\cdot)\nabla u(\cdot,t)):\nabla\ddot{u}(\cdot,t)\,dx,$$
$$J_3=\int_\Omega\big(\int\limits_{0}^{t}F(\cdot,t-\tau)\nabla u(\cdot,\tau)d\tau\big):\nabla\ddot{u}(\cdot,t)\,dx,$$
$$J_4=-s\int_{\Gamma_N}|\ddot{u}(\cdot,t)|^2\,d\sigma-s\gamma\int_{\Gamma_N}\dot{u}(\cdot,t)\ddot{u}(\cdot,t)\,d\sigma.$$
Having in mind (\ref{3.9}) and the following identities
\begin{equation}\label{3.10}
\begin{aligned}
&\frac{d}{dt}\int_\Omega(G(\cdot,0)\nabla u(\cdot,t)):\nabla\dot{u}(\cdot,t)\,dx\\
=&\int_\Omega(G(\cdot,0)\nabla u(\cdot,t)):\nabla\ddot{u}(\cdot,t)\,dx\,+\int_\Omega(G(\cdot,0)\nabla \dot{u}(\cdot,t)):\nabla\dot{u}(\cdot,t)\,dx
\end{aligned}
\end{equation}
\begin{equation}\label{3.11}
\begin{aligned}
&\frac{d}{dt}\int_\Omega(C(\cdot)\nabla u(\cdot,t)) :\nabla\dot{u}(\cdot,t)\,dx\\
=&\gamma\int_\Omega(C(\cdot)\nabla u(\cdot,t)):\nabla\ddot{u}(\cdot,t)\,dx\,+\gamma\int_\Omega(C(\cdot)\nabla \dot{u}(\cdot,t)):\nabla\dot{u}(\cdot,t)\,dx
\end{aligned}
\end{equation}
\begin{equation}\label{3.12}
\begin{aligned}
&\frac{d}{dt}\int_\Omega\Big(\int\limits_{0}^{t}F(\cdot,t-\tau)\nabla u(\cdot,\tau)\,d\tau\Big):\nabla\dot{u}(\cdot,t)\,dx\\
=&\int_\Omega\Big(\int\limits_{0}^{t}F(\cdot,t-\tau)\nabla u(\cdot,\tau)\,d\tau\Big):\nabla\ddot{u}(\cdot,t)\,dx\\
&+\int_\Omega\Big(\int\limits_{0}^{t}\dot{F}(\cdot,t-\tau)\nabla u(\cdot,\tau)\,d\tau\Big):\nabla\dot{u}(\cdot,t)\,dx+\int_\Omega\Big(F(\cdot,0)\nabla u(\cdot,t)\Big):\nabla\dot{u}(\cdot,t)\,dx,
\end{aligned}
\end{equation}
we obtain
\begin{equation}\label{3.13}
\begin{aligned}
\frac{d}{dt}K(t,u)=&-\gamma\int_\Omega |\ddot{u}(\cdot,t)|^2\,dx\,-\int_\Omega \big(G(\cdot,0)\nabla \dot{u}(\cdot,t)\big):\nabla\ddot{u}(\cdot,t)\,dx\\
&+\gamma\int_\Omega \big(C(\cdot)\nabla \dot{u}(\cdot,t)\big):\nabla\dot{u}(\cdot,t)\,dx\,-\int_\Omega \big(F(\cdot,0)\nabla u(\cdot,t)\big):\nabla\dot{u}(\cdot,t)\,dx\\
&-\int_\Omega\Big(\int\limits_{0}^{t}\dot{F}(\cdot,t-\tau)\nabla u(\cdot,\tau)\,d\tau\Big):\nabla\dot{u}(\cdot,t)\,dx\\
&-s\int_{\Gamma_N}|\ddot{u}(\cdot,t)|^2\,d\sigma-s\gamma\int_{\Gamma_N}\dot{u}(\cdot,t)\ddot{u}(\cdot,t)\,d\sigma.
\end{aligned}
\end{equation}
Now multiplying (\ref{3.5}) by $\dot{u}(\cdot,t)$ and integrating in $\Om$, we have
\begin{equation}\label{3.14}
\begin{aligned}
\frac{d}{dt}\int_\Om \ddot{u}(\cdot,t)\dot{u}(\cdot,t)\,dx=J_5+J_6+J_7,
\end{aligned}
\end{equation}
where
$$J_5=\int_\Om|\ddot{u}(\cdot,t)|^2\,dx-\int_\Om\big(C(\cdot)\nabla\dot{u}(\cdot,t)\big):\nabla\dot{u}(\cdot,t)\,dx,$$
$$J_6=\int_\Om\Big(\int\limits_{0}^{t}\dot{G}(\cdot,t-\tau)\nabla u(\cdot,\tau)\,d\tau\Big):\nabla\dot{u}(\cdot,t)\,dx,$$
$$J_7=\int_\Om\big(G(\cdot,0)\nabla u(\cdot,t)\big):\nabla\dot{u}(\cdot,t)\,dx-s\int_{\Gamma_N}\dot{u}(\cdot,t)\ddot{u}(\cdot,t)\,d\sigma.$$
Further, by the definition of $I(t,u)$, (\ref{3.14}) and Lemma \ref{ident1}, we have
\begin{equation}\label{3.15}
\begin{aligned}
\frac{d}{dt}I(t,u)=&\int_\Om |\ddot{u}(\cdot,t)|^2\,dx\,-\int_\Om\big(C(\cdot)\nabla\dot{u}(\cdot,t)\big):\nabla\dot{u}(\cdot,t)\,dx\\
&-\frac{1}{2}\int_\Om \big(\dot{G}(\cdot,t)\nabla u(\cdot,t)\big):\nabla u(\cdot,t)\,dx\,+\frac{1}{2}\int_\Om \ddot{G}\Box\partial u\,dx\\
&-s\int_{\Gamma_N}\dot{u}(\cdot,t)\ddot{u}(\cdot,t)\,d\sigma.
\end{aligned}
\end{equation}
Finally, putting together (\ref{3.13}) and (\ref{3.15}), the proof is complete.
\eproof
\begin{lemma}\label{ident4}
\begin{equation}\label{3.16}
\begin{aligned}
&\frac{d}{dt}\int_\Om u(\cdot,t)\dot{u}(\cdot,t)\,dx\\
=&\int_\Om |\dot{u}(\cdot,t)|^2\,dx-\int_\Omega\{(C(\cdot)-\int_0^t G(\cdot,\tau)\,d\tau)\nabla u(\cdot,t)\}:\nabla u(\cdot,t)\,dx\\
&-\int_\Om\{\int_0^t G(\cdot,t-\tau)(\nabla u(\cdot,t)-\nabla u(\cdot,\tau))d\tau\}:\nabla u(\cdot,t)\, dx-s\int_{\Gamma_N}\dot{u}(\cdot,t)u(\cdot,t)\,d\sigma.
\end{aligned}
\end{equation}
\end{lemma}
\pf. Multiplying (\ref{mixed prob}) by $u(\cdot,t)$ and integrating in $\Om$, we have \eqref{3.16}.
\eproof


\subsection{Basic inequalities}\label{basic ineq}${}$
\noindent
In order to estimate $\frac{d}{dt}\mathcal{L}(t,u)$ from above, we proceed as follows. We estimate $$
\frac{d}{dt} N_1 [ B(t,u)+2\hat{c}_\delta(1+t)^{-p}E(0,u)],\,\,
\frac{d}{dt}\int_\Omega \omega\dot{u}(\cdot,t)u(\cdot,t)\,dx,\,\,\frac{d}{dt}\big\{K(t,u)+(\gamma-\theta)I(t,u)\big\}
$$
which are given below in \eqref{3.29}, \eqref{3.24} and \eqref{3.28}, respectively.
In this subsection, we will use the following convention for positive constants describing inequalities. Namely, $c_j>0,\,j=1,2,\dots$
are all intermediate constants independent of $u$ which are used in deriving some estimates. Once estimates are derived, we use $\hat c_j>0,\, j=1,2,\dots$ for derived estimates which are independent of $u$. We note that we even use these latter constants in proofs of inequalities if we do not estimate any further.
\medskip

\begin{lemma}\label{ident5}
\begin{equation}\label{3.17}
\begin{aligned}
E(t,u)+E(t,\dot{u})\leq \hat{c}_1E(0),
\end{aligned}
\end{equation}
where  $$E(0):=E(0,u)+E(0,\dot{u}).$$
\end{lemma}
\pf. From \eqref{3.2}, it is clear that $E(t,u)\leq E(0,u)$.
Now, recall
\begin{equation*}
\begin{aligned}
B(t,u):= N_3E(t,u)+E(t,\dot u)-\int_\Omega\big(G(\cdot,t)\nabla u(\cdot,0)):\nabla\dot{u}(\cdot,t)\,dx,
\end{aligned}
\end{equation*}
where $N_3$ is large enough such that
\begin{equation}\label{3.18}
\begin{cases}
\begin{aligned}
&\frac{1}{2}N_3 E(0,u)+\frac{1}{2}E(0,\dot u)\leq B(0,u)\leq 2N_3E(0,u)+2E(0,\dot u)\\
&\frac{1}{2}N_3 E(t,u)+\frac{1}{2}E(t,\dot u)\leq B(t,u)+\hat{c}_\delta(1+t)^{-p}E(0,u).
\end{aligned}
\end{cases}
\end{equation}

From \eqref{3.3}, \eqref{1.3} and \eqref{3.18}, we have that
\begin{equation}\label{3.19}
\begin{aligned}
&\int_\Omega\big(G(\cdot,t)\nabla u(\cdot,0)):\nabla\ddot{u}(\cdot,t)\,dx-\frac{1}{2}\int_\Omega\big(G(\cdot,t)\nabla \dot{u}(\cdot,t)\big):\nabla \dot{u}(\cdot,t)\,dx\\
&-\frac{d}{dt} \int_\Omega\big(G(\cdot,t)\nabla u(\cdot,0)):\nabla\dot{u}(\cdot,t)\,dx\\
\leq&|\int_\Omega\big(\dot G(\cdot,t)\nabla u(\cdot,0)):\nabla\dot{u}(\cdot,t)\,dx|-\frac{1}{2}\int_\Omega\big(G(\cdot,t)\nabla \dot{u}(\cdot,t)\big):\nabla \dot{u}(\cdot,t)\,dx\\
\leq & -\frac{1}{4}\alpha_0\int_\Omega|\nabla\dot{u}(\cdot,t)|^2\,dx+c_1\kappa_4(1+t)^{-p}E(0,u)\\
\leq & c_1\kappa_4(1+t)^{-p}E(0,u).
\end{aligned}
\end{equation}
From \eqref{3.19}, it is easy to obtain that
\begin{equation*}
\begin{aligned}
\frac{d}{dt} B(t,u)&\leq c_1\kappa_4(1+t)^{-p}E(0,u).
\end{aligned}
\end{equation*}
Then, integrating
\begin{equation*}
\begin{aligned}
\frac{d}{d\tau} B(\tau,u)&\leq c_1\kappa_4(1+\tau)^{-p}E(0,u)
\end{aligned}
\end{equation*}
over $(0,t)$ and using the first inequality of \eqref{3.18}, we have
\begin{equation*}
\begin{aligned}
B(t,u)\leq B(0,u)+ c_2E(0,u)\leq c_3E(0).
\end{aligned}
\end{equation*}
Combining the second inequality of \eqref{3.18} and this inequality, we derive that
\begin{equation*}
\begin{aligned}
\frac{1}{2} E(t,u)+\frac{1}{2}E(t,\dot u)\leq & \frac{1}{2}N_3 E(t,u)+\frac{1}{2}E(t,\dot u)\\
\leq & B(t,u)+\hat{c}_\delta(1+t)^{-p}E(0,u)\\
\leq & c_3E(0)+\hat{c}_\delta(1+t)^{-p}E(0)\\
\leq & c_4E(0)
\end{aligned}
\end{equation*}
which implies \eqref{3.17}.
\eproof

\begin{lemma}\label{ident10}
\begin{equation}\label{3.29}
\begin{aligned}
&\frac{d}{dt} N_1 [ B(t,u)+2\hat{c}_\delta(1+t)^{-p}E(0,u)]\\
\leq &-\frac{N_1\kappa_2}{2}\int_\Omega G^{1+\frac{1}{p}}\Box\partial u(\cdot,t)\,dx-\frac{N_1\kappa_2}{2}\int_\Omega G^{1+\frac{1}{p}}\Box\partial \dot{u}(\cdot,t)\,dx\\
&-N_1s\int_{\Gamma_N}|\dot{u}(\cdot,t)|^2\,d\sigma-N_1s\int_{\Gamma_N}|\ddot{u}(\cdot,t)|^2\,d\sigma+c_0\kappa_4(1+t)^{-p}E(0),
\end{aligned}
\end{equation}
\end{lemma}
\pf. From \eqref{3.2} and \eqref{3.3}, we have
\begin{equation}\label{r.1}
\begin{aligned}
\frac{d}{dt} B(t,u)
=& N_3\frac{d}{dt} E(t,u)+\frac{d}{dt} E(t,\dot{u})-\frac{d}{dt} \int_\Omega\big(G(\cdot,t)\nabla u(\cdot,0)):\nabla\dot{u}(\cdot,t)\,dx\\
= &\int_\Omega\big(G(\cdot,t)\nabla u(\cdot,0)):\nabla\ddot{u}(\cdot,t)\,dx-\frac{1}{2}\int_\Omega\big(G(\cdot,t)\nabla \dot{u}(\cdot,t)\big):\nabla \dot{u}(\cdot,t)\,dx\\
&-\frac{d}{dt} \int_\Omega\big(G(\cdot,t)\nabla u(\cdot,0)):\nabla\dot{u}(\cdot,t)\,dx
+\frac12\int_\Omega\dot{G}\Box\partial \dot{u}(\cdot,t)\,dx\\
&-s\int_{\Gamma_N}|\ddot{u}(\cdot,t)|^2\,d\sigma-\frac{N_3}{2}\int_\Omega\big(G(\cdot,t)\nabla u(\cdot,t)\big):\nabla u(\cdot,t)\,dx\\
&+\frac{N_3}{2}\int_\Omega\dot{G}\Box\partial u(\cdot,t)\,dx-N_3s\int_{\Gamma_N}|\dot{u}(\cdot,t)|^2\,d\sigma.
\end{aligned}
\end{equation}
Combining this with \eqref{1.2} and \eqref{3.19}, we can easily have \eqref{3.29}.\eproof

For further estimates, we cite the following lemmas given as Lemma 2.2 and Lemma 2.6 in \cite{R-L}.
\begin{lemma}\label{ident6}
For any $u,v\in C^1([0,t];H^1(\Omega))$ we have
\begin{equation}\label{3.20}
\begin{aligned}
&|\int_\Omega\big( \int\limits_{0}^{t}G(\cdot, t-\tau)(\nabla u(\cdot,t)-\nabla u(\cdot,\tau)) d\tau \big) :\nabla v(\cdot,t)\, dx|\\
\leq &
\hat{c}_3\big(\int_\Omega (G^{1+\frac{1}{p}}\Box\partial u)(\cdot,t)\,dx\big)^\frac{1}{2}\big(\int_\Omega|\nabla v(\cdot,t)|^2\,dx  \big)^\frac{1}{2}.
\end{aligned}
\end{equation}
\end{lemma}
We use Lemma 2.6. in \cite{R-L} and \eqref{3.17} to obtain the following lemma.
\begin{lemma}\label{ident7}
\begin{equation}\label{3.21}
\begin{aligned}
\int_\Omega (G\Box\partial u)(\cdot,t)\,dx &\leq
\big(\int_\Omega (G^{1-r}\Box\partial u)(\cdot,t)\,dx\big)^{\frac{1}{2^m}}\big(\int_\Omega (G^{1+\frac{1}{p}}\Box\partial u)(\cdot,t)\,dx\big)^{1-\frac{1}{2^m}}\\
&\leq \hat{c}_4\big(\int_\Omega (G^{1+\frac{1}{p}}\Box\partial u)(\cdot,t)\,dx\big)^{1-\frac{1}{2^m}},
\end{aligned}
\end{equation}
where $r=\frac{2^m-1}{p}<1$.
\end{lemma}
Now, from \eqref{3.20}, we have that
\begin{equation}\label{3.22}
\begin{aligned}
&|\int_\Omega\big( \int\limits_{0}^{t}G(\cdot, t-\tau)(\nabla u(\cdot,t)-\nabla u(\cdot,\tau)) d\tau \big) :\nabla u(\cdot,t)\, dx|\\
\leq &
\frac{8\hat{c}_6}{\mu_0}\int_\Omega (G^{1+\frac{1}{p}}\Box\partial u)(\cdot,t)\,dx+\frac{\mu_0}{8}\int_\Omega|\nabla u(\cdot,t)|^2\,dx.
\end{aligned}
\end{equation}
By the trace Theorem and Poincar\'e inequality, we obtain that
\begin{equation}\label{3.23}
\begin{aligned}
\int_{\Gamma_N}s\dot{u}(\cdot,t)u(\cdot,t)\,d\sigma
\leq \hat{c}_7s^2\int_{\Gamma_N}|\dot{u}(\cdot,t)|^2\,d\sigma+\frac{\mu_0}{8}\int_\Omega|\nabla u(\cdot,t)|^2\,dx.
\end{aligned}
\end{equation}

Combining \eqref{3.16}, \eqref{3.22} and \eqref{3.23}, we have the following.
\begin{lemma}\label{ident8}
\begin{equation}\label{3.24}
\begin{aligned}
\frac{d}{dt}\int_\Omega \omega\dot{u}(\cdot,t)u(\cdot,t)\,dx
\leq &\,\hat{c}_5\omega\int_\Omega|\nabla\dot{u}(\cdot,t)|^2\,dx-\frac{6\mu_0\omega}{8}\,\int_\Omega|\nabla u(\cdot,t)|^2\,dx\\
&+\frac{8\hat{c}_6\omega}{\mu_0}\int_\Omega (G^{1+\frac{1}{p}}\Box\partial u)(\cdot,t)\,dx\,+\hat{c}_7\omega s^2\int_{\Gamma_N}|\dot{u}(\cdot,t)|^2\,d\sigma,
\end{aligned}
\end{equation}
where $\hat{c}_5$ is the Poincar\'e inequality constant.
\end{lemma}
Now, let us estimate some terms in \eqref{3.7}, we use \eqref{1.2}, \eqref{1.3} and \eqref{3.17} to obtain that
\begin{equation}\label{3.25}
\begin{aligned}
-\frac{1}{2}(\gamma-\theta)\int_\Omega\big(\dot{G}(\cdot,t)\nabla u(\cdot,t)\big):\nabla u(\cdot,t)\,dx
\leq & c_5\int_\Omega\big(G(\cdot,t)\nabla u(\cdot,t)\big):\nabla u(\cdot,t)\,dx\\
\leq & \hat{c}_8(1+t)^{-p}.
\end{aligned}
\end{equation}
From \eqref{1.2}, \eqref{1.3} and \eqref{3.17}, we have that
\begin{equation}\label{3.26}
\begin{aligned}
&\frac12(\gamma-\theta)\int_\Omega\ddot{G}\Box\partial u(\cdot,t)\,dx
-\int_\Omega\big(F(\cdot,t)\nabla u(\cdot,t)\big):\nabla\dot{u}(\cdot,t)\,dx\\
\leq &\hat{c}_9\int_\Omega G^{1+\frac{1}{p}}\Box\partial u(\cdot,t)\,dx\\
&+\hat{c}_9\big(\int_\Omega\big(G(\cdot,t)\nabla u(\cdot,t)\big):\nabla u(\cdot,t)\,dx\big)^{\frac{1}{2}}
\big(\int_\Omega\big(G(\cdot,t)\nabla \dot{u}(\cdot,t)\big)\,dx\big)^{\frac{1}{2}}\\
\leq & \hat{c}_9\int_\Omega G^{1+\frac{1}{p}}\Box\partial u(\cdot,t)\,dx+\hat{c}_{10}(1+t)^{-p}.
\end{aligned}
\end{equation}
By \eqref{1.2} and \eqref{3.20}, we obtain that
\begin{equation}\label{3.27}
\begin{aligned}
&\int_\Omega\big\{\int_0^t\dot{F}(\cdot,t-\tau)\nabla\big(u(\cdot,t)-u(\cdot,\tau)\big)\,d\tau\big\}:\nabla\dot{u}(\cdot,t)\,dx\\
\leq &\hat{c}_{11}\big(\int_\Omega G^{1+\frac{1}{p}}\Box\partial u(\cdot,t)\,dx\big)^\frac{1}{2}\big(\int_\Omega|\nabla \dot{u}(\cdot,t)|^2\,dx  \big)^\frac{1}{2}\\
\leq &\frac{4\hat{c}_{11}^2}{\alpha_0}
\int_\Omega G^{1+\frac{1}{p}}\Box\partial u(\cdot,t)\,dx+\frac{\alpha_0}{4}\int_\Omega|\nabla \dot{u}(\cdot,t)|^2\,dx.
\end{aligned}
\end{equation}

Combining \eqref{3.25}, \eqref{3.26}, \eqref{3.27} and using Young's inequality, we have the following.
\begin{lemma}\label{ident9}
\begin{equation}\label{3.28}
\begin{aligned}
&\frac{d}{dt}\big\{K(t,u)+(\gamma-\theta)I(t,u)\big\}\\
\leq&-\theta\int_\Omega|\ddot{u}(\cdot,t)|^2\,dx+(\theta\beta_0-\frac{3\alpha_0}{4})\int_\Omega |\nabla\dot{u}(\cdot,t)|^2\,dx\\
&+(\hat{c}_9+\frac{4\hat{c}_{11}^2}{\alpha_0})\int_\Omega G^{1+\frac{1}{p}}\Box\partial u(\cdot,t)\,dx+(\hat{c}_8+\hat{c}_{10})(1+t)^{-p}\\
&+s(\theta+2\gamma)^2\int_{\Gamma_N}|\dot{u}(\cdot,t)|^2\,d\sigma.
\end{aligned}
\end{equation}
\end{lemma}

\section{Estimates for $\mathcal{L}(t,u)$}${}$
\noindent
In this section, we use a convention for positive constants describing inequalities similar to the previous one in Subsection \ref{basic ineq}. The only difference is that we use $\tilde c_j>0,\,j=1,2,\dots$ instead of $\hat c_j>0,\,j=1,2,\dots$.

\subsection{Upper and lower bounds}

In this subsection we bound $\mathcal{L}(t,u)$ from below and above. We will first estimate $N_3E(t,u)+E(t,\dot{u})$ in \eqref{4.0} and then bound $B(t,u)+2\hat{c}_\delta(1+t)^{-p}E(0,u)$ in \eqref{4.1}.
From \eqref{1.5}, we have the following estimate.
\begin{lemma}\label{lem4.1}
There exists two positive constants $\tilde{c}_1< \tilde{c}_2$ such that
\begin{equation}\label{4.0}
\begin{aligned}
\tilde{c}_1  M(t,u)\leq N_3 E(t,u)+E(t,\dot{u})\leq \tilde{c}_2  M(t,u),
\end{aligned}
\end{equation}
where
\begin{equation*}
\begin{aligned}
M(t,u)=R(t,u)+\int_\Omega G\Box\partial u(\cdot,t)\,dx+\int_\Omega G\Box\partial \dot{u}(\cdot,t)\,dx
\end{aligned}
\end{equation*}
and
\begin{equation*}
\begin{aligned}
R(t,u)=\int_\Omega|\ddot{u}(\cdot,t)|^2\,dx+\int_\Omega |\nabla u(\cdot,t)|^2\,dx+\int_\Omega |\nabla\dot{u}(\cdot,t)|^2\,dx.
\end{aligned}
\end{equation*}
\end{lemma}

Also, \eqref{B} gives the following estimate.
\begin{lemma}\label{lem4.2}
There exists two positive constants $\tilde{c}_3< \tilde{c}_4$ such that
\begin{equation}\label{4.1}
\begin{aligned}
\tilde{c}_3 \tilde{M}(t,u)\leq  B(t,u)+2\hat{c}_\delta(1+t)^{-p}E(0,u) \leq \tilde{c}_4 \tilde{M}(t,u),
\end{aligned}
\end{equation}
where $\tilde{M}(t,u)=M(t,u)+(1+t)^{-p}E(0,u)$.
\end{lemma}

Now to obtain the upper bound of $|K(t,u)|$ and $|I(t,u)|$, we use \eqref{1.2} and \eqref{3.20} to estimate the following term.
\begin{equation}\label{4.2}
\begin{aligned}
&\big|\int_\Om\big(\int\limits_{0}^{t}F(\cdot,t-\tau)\nabla u(\cdot,\tau)\,d\tau \big):\nabla\dot{u}(\cdot,t)\,dx\big|\\
\leq &\big|\int_\Om\big(\int\limits_{0}^{t}(\gamma G-\dot{G})(\cdot,t-\tau)\Big(\nabla u(\cdot,t)\,-\nabla u(\cdot,\tau)\Big)\,d\tau \big):\nabla\dot{u}(\cdot,t)\,dx\big|\\
&+\big|\int_\Om\big(\int\limits_{0}^{t}(\gamma G-\dot{G})(\cdot,t-\tau)d\tau\,\nabla u(\cdot,t) \big):\nabla\dot{u}(\cdot,t)\,dx\big|\\
\leq &(\gamma + \kappa_1)\big|\int_\Om\big(\int\limits_{0}^{t}G(\cdot,t-\tau)\Big(\nabla u(\cdot,t)\,-\nabla u(\cdot,\tau)\Big)\,d\tau \big):\nabla\dot{u}(\cdot,t)\,dx\big|\\
&+(\gamma + \kappa_1)\big|\int_\Om\big(\int\limits_{0}^{t}G(\cdot,t-\tau)d\tau\,\nabla u(\cdot,t) \big):\nabla\dot{u}(\cdot,t)\,dx\big|\\
\leq &c_6\gamma\int_\Omega G\Box\partial u(\cdot,t)\,dx+c_6\gamma\int_\Omega |\nabla u(\cdot,t)|^2\,dx+c_6\gamma\int_\Omega |\nabla\dot{u}(\cdot,t)|^2\,dx.
\end{aligned}
\end{equation}
Combining \eqref{4.0}, \eqref{4.1} and \eqref{4.2}, simple computations give the following estimate of $ \mathcal{L}(t,u) $.

\begin{lemma}\label{lem4.3}
There exists two positive constants $\tilde{c}_5< \tilde{c}_6$ such that
\begin{equation}\label{4.3}
\begin{aligned}
\tilde{c}_5 \tilde{M}(t,u)\leq \mathcal{L}(t,u) \leq \tilde{c}_6\tilde{M}(t,u),
\end{aligned}
\end{equation}
if we take $N_1$ large enough.
\end{lemma}

\subsection{Estimating the time derivative from above} \label{bound derivative}

To bound $\frac{d}{dt}\mathcal{L}(t,u)$ from above, we take $\theta=\frac{\alpha_0}{4\beta_0}, \gamma=2\theta, \omega=\frac{\alpha_0}{4\hat{c}_5}$.
Combining \eqref{3.24}, \eqref{3.28} and \eqref{3.29}, we obtain for $N_1$ large enough that

\begin{equation}\label{4.4}
\begin{aligned}
\frac{d}{dt}\mathcal{L}(t,u)\leq &(\hat{c}_8+\hat{c}_{10}+\hat{c}_2E(0))(1+t)^{-p}-\frac{\alpha_0}{4\beta_0}\int_\Om |\ddot{u}(\cdot,t)|^2\,dx
-\frac{3\mu_0\alpha_0}{16}\,\int_\Omega|\nabla u(\cdot,t)|^2\,dx\\
&-\frac{\alpha_0}{4}\int_\Om |\nabla\dot{u}(\cdot,t)|^2\,dx\,
-\big(\frac{N_1\kappa_2}{2}-\frac{2\hat{c}_6\alpha_0}{\hat{c}_5\mu_0}-\hat{c}_9-\frac{4\hat{c}_{11}^2}{\alpha_0}\big)\int_\Omega G^{1+\frac{1}{p}}\Box\partial u(\cdot,t)\,dx\\
&-\frac{N_1\kappa_2}{2}\int_\Omega G^{1+\frac{1}{p}}\Box\partial \dot{u}(\cdot,t)\,dx
-s(N_1-25\theta^2-c_3s\omega)\int_{\Gamma_N}|\dot{u}(\cdot,t)|^2\,d\sigma\\
\leq & c_7(1+t)^{-p}-c_8H(t,u),
\end{aligned}
\end{equation}
where
\begin{equation*}
\begin{aligned}
H(t,u)=(1+t)^{-p}+R(t,u)+\int_\Omega G^{1+\frac{1}{p}}\Box\partial u(\cdot,t)\,dx+\int_\Omega G^{1+\frac{1}{p}}\Box\partial \dot{u}(\cdot,t)\,dx.
\end{aligned}
\end{equation*}

Let us denote $\delta_m=\frac{1}{2^m-1}$. From \eqref{4.3}, \eqref{3.17} and \eqref{3.21}, we have
\begin{equation}\label{4.5}
\begin{aligned}
\mathcal{L}(t,u)^{(1+\delta_m)} \leq & \tilde{c}_6^{(1+\delta_m)} \tilde{M}(t,u)^{(1+\delta_m)}\\
\leq & c_9R(t,u)^{\delta_m}R(t,u)+c_9\big(\int_\Omega G\Box\partial u(\cdot,t)\,dx\big)^{(1+\delta_m)}\\
&+c_98\big(\int_\Omega G\Box\partial \dot{u}(\cdot,t)\,dx\big)^{(1+\delta_m)}+c_{9}(1+t)^{-p(1+\delta_m)}\\
\leq & c_{10}R(t,u)+c_{10}\int_\Omega G^{1+\frac{1}{p}}\Box\partial u(\cdot,t)\,dx+c_{10}\int_\Omega G^{1+\frac{1}{p}}\Box\partial \dot{u}(\cdot,t)\,dx\\
&+c_{10}(1+t)^{-p}\\
=& c_{10} H(t,u).
\end{aligned}
\end{equation}

\medskip
Combining \eqref{4.4} and \eqref{4.5}, we have the following lemma.
\begin{lemma}\label{lem4.4}
There exists two positive constants $\tilde{c}_7$, $\tilde{c}_8$ such that
\begin{equation}\label{4.6}
\begin{aligned}
\frac{d}{dt}\mathcal{L}(t,u)\leq &\tilde{c}_7(1+t)^{-p}-\tilde{c}_{8}\mathcal{L}(t,u)^{(1+\delta_m)}\\
\leq &\tilde{c}_7(1+t)^{-p_m-1}-\tilde{c}_{8}\mathcal{L}(t,u)^{(1+\frac{1}{p_m})},
\end{aligned}
\end{equation}
where $m$ is a positive integer such that $p_m=2^m-1<p-1$.
\end{lemma}

\section{Proof of Theorem \ref{mainresult2}}${}$
In this section we give a proof of Theorem \ref{mainresult2}. Since the proof is not only almost the same as that of Theorem \ref{mainresult} but also simpler than that of Theorem \ref{mainresult}, we only give necessary changes which we need to make in the proof of Theorem \ref{mainresult} to give a proof of Theorem \ref{mainresult2}.

Note that due to \eqref{derivatives} and \eqref{exp decay of the |G| and |dot G|},
$G(t)$ does not satisfy \eqref{1.2} and \eqref{1.3}.
Nevertheless, by using \eqref{derivatives} and \eqref{exp decay of the |G| and |dot G|}, we can obtain the uniform decay rate of exponential order for the solution of (MP)'.
Now, we define
\begin{equation*}
\begin{aligned}
\tilde{B}(t,u):= N_4E(t,u)+E(t,\dot u)-\int_\Omega\big(G(\cdot,t)\nabla u(\cdot,0)):\nabla\dot{u}(\cdot,t)\,dx,
\end{aligned}
\end{equation*}
where $N_4$ is large enough such that
\begin{equation}\label{Be}
\begin{cases}
\begin{aligned}
&\frac{1}{2}N_4 E(0,u)+\frac{1}{2}E(0,\dot u)\leq \tilde{B}(0,u)\leq 2N_4E(0,u)+2E(0,\dot u)\\
&\frac{1}{2}N_4 E(t,u)+\frac{1}{2}E(t,\dot u)\leq \tilde{B}(t,u)+\hat{c}_\delta e^{-\kappa_4t}E(0,u).
\end{aligned}
\end{cases}
\end{equation}
We also let
\begin{equation}\label{mathcal Le}
\mathcal{L}_e(t,u):=N_1 [ \tilde{B}(t,u)+2\hat{c}_\delta e^{-\kappa_4t}E(0,u)]+K(t,u)+(\gamma-\theta)I(t,u)+\omega\int_\Omega\dot{u}u\,dx.
\end{equation}

By replacing \eqref{4.6}, we will prove
\begin{equation}\label{5.11}
\frac{d}{dt}\mathcal{L}_e(t,u)\leq e^{-\tilde{\kappa}_4 t}-b_2\mathcal{L}_e(t,u),\,\,t>0,
\end{equation}
where $b_2>0$ is a constant which does not depend on $u$. The lower bound of $\mathcal{L}_e(t,u)$ similar to \eqref{4.3} can be obtained easier than before.

The details of the proof of \eqref{5.11} are similar to that of Theorem 1.2 of \cite{NO}.
Thus, from \eqref{4.3} and \eqref{5.11}, we have
\begin{equation}\label{estimate of U(T) exp_decay}
\mathcal{L}_e(t,u)\le \mathcal{L}_e(0,u)\, a_3e^{-b_3 t}
\end{equation}
with some constants $a_3, b_3>0$ independent of $t$.

In the following, we state the differences of the estimates using \eqref{derivatives}, \eqref{exp decay of the |G| and |dot G|} and the previous estimates using \eqref{1.2}, \eqref{1.3}. The new constants $c_j$'s and $\kappa_j$'s in the forthcoming estimates are all positive constants independent of $u$.

First, \eqref{3.20} changes to
\begin{equation}\label{5.21}
\begin{aligned}
&|\int_\Omega\big( \int\limits_{0}^{t}G(\cdot, t-\tau)(\nabla u(\cdot,t)-\nabla u(\cdot,\tau)) d\tau \big) :\nabla v(\cdot,t)\, dx|\\
\leq &
\hat{c}_3\big(\int_\Omega (G\Box\partial u)(\cdot,t)\,dx\big)^\frac{1}{2}\big(\int_\Omega|\nabla v(\cdot,t)|^2\,dx  \big)^\frac{1}{2}.
\end{aligned}
\end{equation}
\eqref{3.29}  changes to
\begin{equation}\label{5.27}
\begin{aligned}
&\frac{d}{dt} N_1 [ \tilde{B}(t,u)+2\hat{c}_\delta e^{-\kappa_4t}E(0,u)]\\
\leq &-\frac{N_1\kappa_2}{2}\int_\Omega G\Box\partial u(\cdot,t)\,dx-\frac{N_1\kappa_2}{2}\int_\Omega G\Box\partial \dot{u}(\cdot,t)\,dx\\
&-N_1s\int_{\Gamma_N}|\dot{u}(\cdot,t)|^2\,d\sigma-N_1s\int_{\Gamma_N}|\ddot{u}(\cdot,t)|^2\,d\sigma+\hat{c}_2\kappa_4e^{-\tilde{\kappa}_4 t}E(0).
\end{aligned}
\end{equation}
\eqref{3.24}  changes to
\begin{equation}\label{5.22}
\begin{aligned}
\frac{d}{dt}\int_\Omega \omega\dot{u}(\cdot,t)u(\cdot,t)\,dx
\leq &\,\hat{c}_5\omega\int_\Omega|\nabla\dot{u}(\cdot,t)|^2\,dx-\frac{6\mu_0\omega}{8}\,\int_\Omega|\nabla u(\cdot,t)|^2\,dx\\
&+\frac{8\hat{c}_6\omega}{\mu_0}\int_\Omega ( G\Box\partial u)(\cdot,t)\,dx\,+\hat{c}_7\omega s^2\int_{\Gamma_N}|\dot{u}(\cdot,t)|^2\,d\sigma.
\end{aligned}
\end{equation}
\eqref{3.25}  changes to
\begin{equation}\label{5.23}
\begin{aligned}
-\frac{1}{2}(\gamma-\theta)\int_\Omega\big(\dot{G}(\cdot,t)\nabla u(\cdot,t)\big):\nabla u(\cdot,t)\,dx
\leq & c_5\int_\Omega\big(G(\cdot,t)\nabla u(\cdot,t)\big):\nabla u(\cdot,t)\,dx\\
\leq & \hat{c}_8e^{-\tilde{\kappa}_4 t}.
\end{aligned}
\end{equation}
\eqref{3.26}  changes to
\begin{equation}\label{5.24}
\begin{aligned}
&\frac12(\gamma-\theta)\int_\Omega\ddot{G}\Box\partial u(\cdot,t)\,dx
-\int_\Omega\big(F(\cdot,t)\nabla u(\cdot,t)\big):\nabla\dot{u}(\cdot,t)\,dx\\
\leq & \hat{c}_9\int_\Omega G^{1+\frac{1}{p}}\Box\partial u(\cdot,t)\,dx+\hat{c}_{10}e^{-\tilde{\kappa}_4 t}.
\end{aligned}
\end{equation}
\eqref{3.27}  changes to
\begin{equation}\label{5.25}
\begin{aligned}
&\int_\Omega\big\{\int_0^t\dot{F}(\cdot,t-\tau)\nabla\big(u(\cdot,t)-u(\cdot,\tau)\big)\,d\tau\big\}:\nabla\dot{u}(\cdot,t)\,dx\\
\leq &\frac{4\hat{c}_{11}^2}{\alpha_0}\int_\Omega G\Box\partial u(\cdot,t)\,dx+\frac{\alpha_0}{4}\int_\Omega|\nabla \dot{u}(\cdot,t)|^2\,dx.
\end{aligned}
\end{equation}
\eqref{3.28}  changes to
\begin{equation}\label{5.26}
\begin{aligned}
&\frac{d}{dt}\big\{K(t,u)+(\gamma-\theta)I(t,u)\big\}\\
\leq&-\theta\int_\Omega|\ddot{u}(\cdot,t)|^2\,dx+(\theta\beta_0-\frac{3\alpha_0}{4})\int_\Omega |\nabla\dot{u}(\cdot,t)|^2\,dx\\
&+(\hat{c}_9+\frac{4\hat{c}_{11}^2}{\alpha_0})\int_\Omega  G\Box\partial u(\cdot,t)\,dx+(\hat{c}_8+\hat{c}_{10})e^{-\tilde{\kappa}_4 t}\\
&+s(\theta+2\gamma)^2\int_{\Gamma_N}|\dot{u}(\cdot,t)|^2\,d\sigma.
\end{aligned}
\end{equation}

Secondly, combining \eqref{5.22}, \eqref{5.26} and \eqref{5.27}, \eqref{4.4}  changes to
\begin{equation}\label{5.28}
\begin{aligned}
&\frac{d}{dt}\mathcal{L}_e(t,u)\\
\leq &(\hat{c}_8+\hat{c}_{10}+\hat{c}_2E(0))e^{-\tilde{\kappa}_4 t}-\frac{\alpha_0}{4\beta_0}\int_\Om |\ddot{u}(\cdot,t)|^2\,dx
-\frac{3\mu_0\alpha_0}{16}\,\int_\Omega|\nabla u(\cdot,t)|^2\,dx\\
&-\frac{\alpha_0}{4}\int_\Om |\nabla\dot{u}(\cdot,t)|^2\,dx\,
-\big(\frac{N_1\kappa_2}{2}-\frac{2\hat{c}_6\alpha_0}{\hat{c}_5\mu_0}-\hat{c}_9-\frac{4\hat{c}_{11}^2}{\alpha_0}\big)\int_\Omega  G\Box\partial u(\cdot,t)\,dx\\
&-\frac{N_1\kappa_2}{2}\int_\Omega G\Box\partial \dot{u}(\cdot,t)\,dx
-s(N_1-25\theta^2-c_3s\omega)\int_{\Gamma_N}|\dot{u}(\cdot,t)|^2\,d\sigma\\
\leq &c_7e^{-\tilde{\kappa}_4 t}-c_8[M(t,u)+e^{-\tilde{\kappa}_4 t}].
\end{aligned}
\end{equation}
\eqref{3.29}  changes to
\begin{equation}\label{5.29}
\begin{aligned}
\tilde{c}_5[M(t,u)+e^{-\tilde{\kappa}_4 t}]\leq \mathcal{L}_e(t,u) \leq \tilde{c}_6[M(t,u)+e^{-\tilde{\kappa}_4 t}].
\end{aligned}
\end{equation}
Combining \eqref{5.28} and \eqref{5.29}, we can derive \eqref{5.11}.

\section {Conclusion}${}$
    We clarified and modified the method of \cite{R-L}, and pushed it further, showing the followings in a unified way. We gave some sufficient conditions for solutions of the initial boundary value problem with mixed type boundary condition for the VID system to have the polynomial decaying property and the exponential decaying property. By considering well known viscoelastic models in mechanics such as the Maxwell model, the SLS model, the Burgers  and Maxwell models and their extended versions, we showed that these provide the VID systems under some condition at $t=0$. Further, we examined whether they satisfy Assumption \ref{exponential} or not. Except the Burgers  and Maxwell models and their extended versions, the VID systems associated to the other models can satisfy Assumption \ref{exponential}. Nevertheless, by parallelly connecting one spring, the VID systems associated to the resulting models satisfy Assumption \ref{exponential}.

   Concerning more general decay rates besides exponential or polynomial rates, see \cite{Massaoudi} and references there in. 

\appendix

\section{Examples of $G$}

\subsection{Separable $G$ with polynomial order decaying factor} \label{secA}${}$
\newline
\indent
First of all, it is not easy to find a separable $G$ coming from a spring-dashpot models satisfying \eqref{1.2} and \eqref{1.3}. To find such a $G(x,t)=G(t) \hat{G }(x)$ with $\hat G(x)$ satisfying the major symmetry and the strong convexity conditions, we look for $0\leq G(t)\in C^2([0,\infty))$ which satisfies the following estimates
\begin{equation}\label{a.1}
-\kappa_1G(t)\le\dot{G}(t)\le-\kappa_2G^{1+\frac{1}{p}}(t),\quad \ddot{G}(t)\le \kappa_3G^{1+\frac{1}{p}}(t)
\end{equation}
which corresponds to \eqref{1.2}. As we have already discussed in Remark \ref{several remarks}, we can assume $G(t)>0,\,t\ge0$.

From $\dot{G}(t)\le-\kappa_2G^{1+\frac{1}{p}}(t)$, we have
\begin{equation}\label{a.2}
\begin{aligned}
p\frac{d}{dt} G(t)^{-\frac{1}{p}}=-\dot GG^{-1-\frac{1}{p}}\geq \kappa_2.
\end{aligned}
\end{equation}
Then, integrating
\begin{equation*}
\begin{aligned}
\frac{d}{dt} G(t)^{-\frac{1}{p}}&\geq \kappa_2/p
\end{aligned}
\end{equation*}
over $[0,t]$, we have
\begin{equation*}
\begin{aligned}
G(t)^{-\frac{1}{p}}\geq G(0)^{-\frac{1}{p}}+\frac{\kappa_2}{p}t.
\end{aligned}
\end{equation*}
Thus, we obtain that
\begin{equation}\label{a.3}
\begin{aligned}
G(t)\leq (G(0)^{-\frac{1}{p}}+\frac{\kappa_2}{p}t)^{-p}.
\end{aligned}
\end{equation}
Combining this with $-\kappa_1G(t)\le\dot{G}(t)$, we have that
\begin{equation}\label{a.4}
\begin{aligned}
|\dot G|+G(t)\leq (1+\kappa_1)(G(0)^{-\frac{1}{p}}+\frac{\kappa_2}{p}t)^{-p}.
\end{aligned}
\end{equation}
Here, observe that if $\kappa_2=p\cdot \hat{\kappa}$ for some $\hat{\kappa}>0$ and $G(0)\leq 1$, then this gives \eqref{1.3}.

From $-\kappa_1G(x,t)\le\dot{G}(x,t)$, we have that
\begin{equation}\label{a.5}
\begin{aligned}
\frac{d}{dt} \ln G(t)=\dot GG^{-1}\geq -\kappa_1.
\end{aligned}
\end{equation}
Then, integrating \eqref{a.5} over $[0,t]$, we have
\begin{equation}\label{a.6}
\begin{aligned}
G(t) &\geq G(0) e^{-\kappa_1t},
\end{aligned}
\end{equation}
which suggest us to look for $G(t)$ as a negative power of a liner funtion of $t$.

Based on this and the above observation, it is natural to look for $G(t)$ in the form
$$ G(t)=(1+at)^{-p},\quad a>0.$$
A direct computation gives that
\begin{equation}\label{a.7}
\begin{aligned}
-apG(t)\leq \dot G(t)=-ap(1+at)^{-p-1}\leq -a p \, G(t)^{1+\frac{1}{p}}
\end{aligned}
\end{equation}
and
\begin{equation}\label{a.8}
\begin{aligned}
\ddot G(t)=ap(p+1)(1+at)^{-p-2}\leq ap(p+1)G(t)^{1+\frac{1}{p}}.
\end{aligned}
\end{equation}
Thus, $G(t)=(1+at)^{-p}$ satisfies \eqref{a.1} and \eqref{a.4}.

\subsection{Spring-dashpot models and exponential decay} \label{secA_M}
The schematic figures of the one space dimensional Maxwell model, standard linear solid (SLS) model and Burgers model are as follow:
\begin{center}
\includegraphics[width=1\textwidth, height=0.25\textheight]{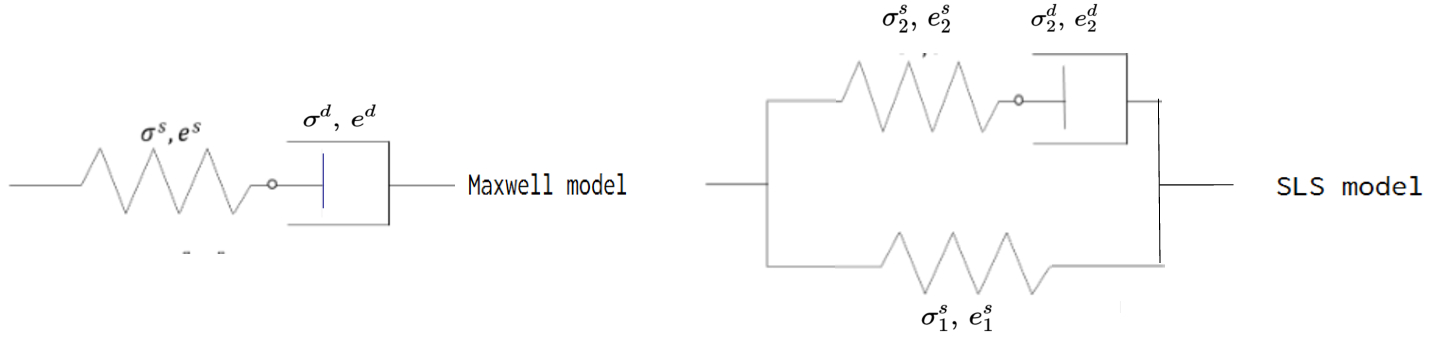}
\end{center}
\begin{center}
\includegraphics[width=0.8\textwidth, height=0.25\textheight]{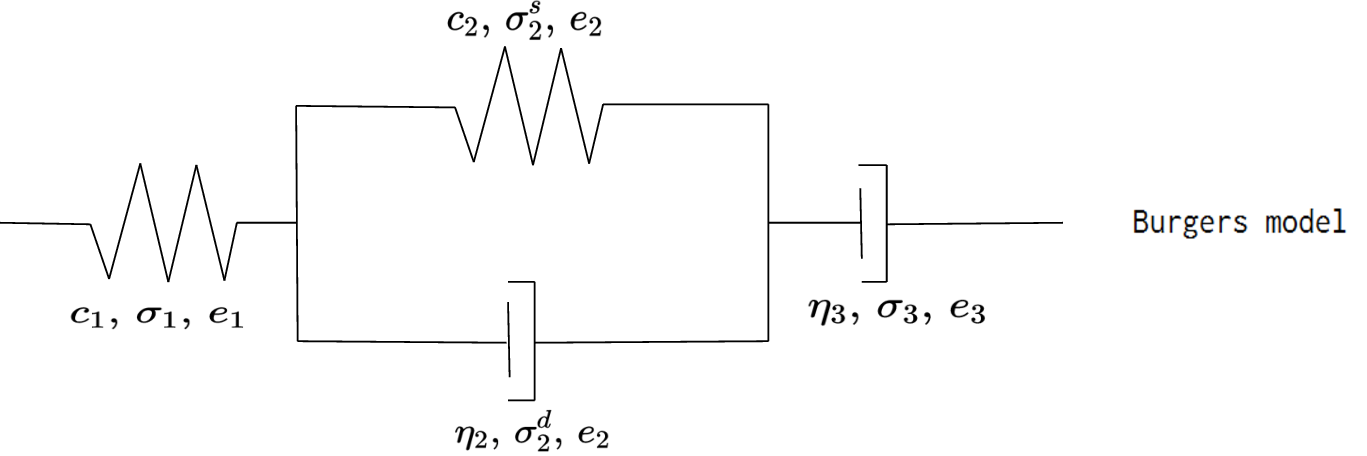}
\end{center}
We will show that these dash-pot models and their extended versions yield VID systems with separable kernels. Also, for the extended Maxwell model and the extended SLS model, their three dimensional generalization also yield VID systems with separable kernels which can satisfy Assumption\ref{exponential}.

\bigskip
\noindent
{\bf Maxwell model:}${}$
\newline
Let $\sigma^s$, $e^s$ and $C_s$ be the stress, strain, and spring constant of the spring, respectively. Also, let
$\sigma^d$, $e^d$ and $\eta$ be the stress, strain and viscosity constant of the dashpot, respectively. Then from Hooke's law, we have
\begin{equation}\label{b.1}
\begin{aligned}
\sigma^s=C_s e^s,
\end{aligned}
\end{equation}
where $\sigma^s$, $e^s$ and $C_s$ are stress, strain, and spring constant of a spring, respectively. Also, for the dashpot, it is usually assumed to have the following relation
\begin{equation}\label{b.2}
\begin{aligned}
\sigma^d=\eta\,\dot e^d.
\end{aligned}
\end{equation}

Since the spring and dashpot are connected in series, the stress $\sigma$ and strain $\epsilon$ of this model are given as
\begin{equation}\label{b.3}
\sigma:=\sigma^s=\sigma^d,\,\,e:=e^s+e^d.
\end{equation}
By Newton's second law, we have
\begin{equation}\label{Newton}
\rho\, \partial_t^2 u=\partial_x \sigma,
\end{equation}
where $\rho$ is the density of the system and $u$ is the displacement of the system which gives the expression $e=\partial_x u$.
Then, from \eqref{b.1} to \eqref{Newton}, we have
the following augmented system
\begin{equation}\label{b.5}
\rho\,\partial_t^2 u = \partial_x(C_s e - C_s e^d ),\,\,
\eta\,\dot{e}^d = C_s e - C_s e^d.
\end{equation}
By assuming $e^d(0)=0$, a simple manipulation yields the constitutive equation
\begin{equation}\label{b.4}
\sigma=C_s e - \int_0^t G(t-\tau)e(\tau)\,
d\tau
\end{equation}
with
\begin{equation}\label{G(t)}
G(t)=a_1 \exp(-b_1 t),
\end{equation}
where $a_1=\eta^{-1}C_s^2,\, b_1=\eta^{-1}C_s$.
If we further assume $e=0$ at $t=0$, then \eqref{b.4} can be written in the form
\begin{equation}\label{eq with relaxation}
\sigma(t)=\int_0^t H(t-\tau)\,\partial_\tau e(\tau)\,d\tau=H(0)e(t)+\int_0^t \dot H(t-\tau) e(\tau)\,d\tau.
\end{equation}
with the relaxation tensor $H(t)=C_s \exp(-b_1 t)$ and
$\dot H(t)=dH/(dt)(t)$.
If the above physcial constants depend on $x$, \eqref{eq with relaxation} becomes
\begin{equation}\label{x dep eq  with relaxation}
\sigma(x,t)=\int_0^t H(x,t-\tau)\,\partial_\tau e(x,\tau)\,d\tau=H(x,0)e(x,t)+\int_0^t \dot H(x,t-\tau)e(x,\tau)\,d\tau.
\end{equation}
Here, note that $C(x)$ and $G(x,t)$ in \eqref{mixed prob} are $C(x)=H(x,0)$ and $G(x,t)=-\dot{H}(x,t)$. Also, $G(t)$ given by \eqref{G(t)} can satisfy Assumption \ref{exponential} in Section \ref{sec1} except (iv)'.

The extended Maxwell model is obtained by connecting the Maxwell model in parallel. Hence, its relaxation tensor is the sum of the relaxation tensor of Maxwell models, and it can satisfy Assumption \ref{exponential} except (iv)'.
\eqref{x dep eq  with relaxation} can be further extended by considering higher space dimensional cases, $n=2, 3$ cases, and also considering anisotropic tensors $C(x)$ and $G(x,t)$. Then, the resulting VID equations can satisfy Assumption \ref{exponential} except (iv)'.  Further, if we connect one spring to this extended Maxwell model, then Assumption
1.2 is satisfied. This is due to the fact that the parallel connection of one spring only adds another stress of the newly connected spring to the original total stress.
Note that for $n=1$ case, the relaxation tensor is just a function. So, it would be more comfortable to called it relaxation kernel. However, even in this case we will call it the relaxation tensor. In this paper, we will use this convention in this paper.

\bigskip
\noindent
{\bf SLS model:}${}$
\newline
Let $\{\sigma_1^s, \sigma_2^s\},\,\{e_1^s,e_2^s\}$ and $\{C_1,C_2\}$ be stresses, strains and spring constants of the springs, respectively. Also, let $\sigma_2^d, e_2^d$ and $\eta_2$ be the stress and strain and viscosity constant of dashpot, respectively. Then, we have
$$
\sigma_1^s=C_1e_1^s,\,\,\sigma_2^s=C_2e_2^s,\,\,\sigma_2^d=\eta_2\dot e_2^d.
$$
Also, the total stress $\sigma$ and the strain $e$ of the model are given as
$$
\sigma:=\sigma_1^s+\sigma_2^d\,\,\text{with}\,\,\sigma_2^s=\sigma_2^d,\,\,e:=e_1^s=e_2^s+e_2^d.
$$
By direct computation as we did for the Maxwell model case, the augmented system is given as
$$
\rho\partial_t^2 u=\partial_x\{(C_1+C_2)e-C_2 e_2^d\}=0,\,\,\eta_2\dot e_2^d=C_2(e-e_2^d),
$$
and by assuming $e_2^d(0)=0$, the constitutive equation is given as
\begin{equation}\label{constitutive SLS}
\sigma(t)=(C_1+C_2)e(t)-\int_0^t\exp\{-b_2(t-s)\}\eta_2^{-1} C_2^2 e(s)\,ds,
\end{equation}
where
$$
a_2:=\eta_2^{-1}C_2^2,\,\,b_2:=\eta_2^{-1} C_2.
$$
Hence $G(x,t)=\exp(b_2t)\eta_2^{-1}C_2^2$ and it can satisfy Assumption \ref{exponential}. Note that $C_1 e(t)$ term is very important for satisfying Assumption 1.2, (iv)'. 

The extended SLS model is obtained by connecting the SLS model in parallel. Hence, the relaxation tensor is the sum of the relaxation tensors of the SLS model, and it can satisfy Assumption \ref{exponential}. 

For the three space dimensional case, we also can say the same thing which we said for the Maxwell model.

\bigskip\noindent
{\bf Burgers model:}${}$
\newline
Let $\{\sigma_1, \sigma_2^s\}$, $\{e_1,e_2\}$ and $\{c_1,c_2\}$ be the stresses, strains and spring constants of the springs, respectively. Also, let
$\{\sigma_2^d,\sigma_3\}$, $\{e_2,e_3\}$ and $\{\eta_2,\eta_3\}$ be the stresses, strains and viscosity constants of the dashpots, respectively.
Note that
\begin{equation}\label{d.1}
\begin{aligned}
\sigma:=\sigma_1=\sigma_2=\sigma_3,
\end{aligned}
\end{equation}
where $\sigma_2=\sigma_2^s+\sigma_2^d$.
By Hooke's law, we have
\begin{equation}\label{d.2}
\sigma_1=c_1e_1,\,\,
\sigma_2^s=c_2e_2.
\end{equation}
Since the stress of a dashpot is proportional to the speed of strain of the dashpot, we have
\begin{equation}\label{d.3}
\sigma_2^d=\eta_2\dot{e_2},\,\,
\sigma_3=\eta_3\dot{e_3}.
\end{equation}

Let $e$ be the total strain given as $e:=e_1+e_2+e_3$. Then, the relation between $\sigma$ and $e$ is
\begin{equation}\label{d.4}
\begin{aligned}
\sigma=\sigma_2=&\sigma_2^s+\sigma_2^d=c_2e_2+\eta_2\dot{e_2}\\
=&c_2(e-e_1-e_3)+\eta_2(\dot{e}-\dot{e_1}-\dot{e_3})\\
=&c_2(e-c_1^{-1}\sigma-e_3)+\eta_2(\dot{e}-c_1^{-1}\dot{\sigma}-\eta_3^{-1}\sigma).
\end{aligned}
\end{equation}
By differentiating  \eqref{d.4} by $t$, we have
\begin{equation}\label{d.5}
\begin{aligned}
\dot{\sigma}=&c_2(\dot{e}-c_1^{-1}\dot{\sigma}-\dot{e_3})+\eta_2(\ddot{e}-c_1^{-1}\ddot{\sigma}-\eta_3^{-1}\dot{\sigma})\\
=&c_2(\dot{e}-c_1^{-1}\dot{\sigma}-\eta_3^{-1}\sigma)+\eta_2(\ddot{e}-c_1^{-1}\ddot{\sigma}-\eta_3^{-1}\dot{\sigma}).
\end{aligned}
\end{equation}

Now, rewrite \eqref{d.5} in the form:
\begin{equation}\label{d.6}
\begin{aligned}
\sigma+B_1\dot{\sigma}+B_2\ddot{\sigma}=\eta_3\dot{e}+\eta_2\eta_3 c_2^{-1}\ddot{e},
\end{aligned}
\end{equation}
where $B_1=\eta_3 c_2^{-1}+\eta_3 c_1^{-1}+\eta_2 c_2^{-1}$ and $B_2=\eta_3 c_1^{-1}\eta_2c_2^{-1}$.
By taking the Laplace transform of \eqref{d.6} with respect to $t$, we have
\begin{equation}\label{d.7}
\begin{aligned}
&\hat{\sigma}(s)+B_1s\hat{\sigma}(s)+B_2s^2\hat{\sigma}(s)-B_1\sigma(0)-B_2\dot{\sigma}(0)-sB_2\sigma(0)\\
=&\eta_3s\hat{e}(s)+\eta_2\eta_3 c_2^{-1}s^2\hat{e}(s)-\eta_3e(0)-\eta_2\eta_3 c_2^{-1}\dot{e}(0)-s\eta_2\eta_3 c_2^{-1}e(0),
\end{aligned}
\end{equation}
where $\hat{\sigma}(s)$ and $\hat{e}(s)$ are the Laplace transforms of $\sigma(t)$ and $e(t)$, respectively.
Here, we assume that the $\sigma(t)$ and $e(t)$ satisfy the following conditions:
\begin{equation}\label{d.8}
B_2\sigma(0)=\eta_2\eta_3 c_2^{-1}e(0),\,\,
B_1\sigma(0)+B_2\dot{\sigma}(0)=\eta_3e(0)+\eta_2\eta_3 c_2^{-1}\dot{e}(0).
\end{equation}
Combining \eqref{d.7} and \eqref{d.8}, we obtain
\begin{equation}\label{d.9}
\begin{aligned}
\hat{\sigma}(s)+B_1s\hat{\sigma}(s)+B_2s^2\hat{\sigma}(s)=\eta_3s\hat{e}(s)+\eta_2\eta_3 c_2^{-1}s^2\hat{e}(s).
\end{aligned}
\end{equation}
Our aim is to find $G(t)$ such that
\begin{equation}\label{d.10n}
\begin{aligned}
\sigma(t)=G(0)e(t)+\int_0^t \dot{G}(t-\tau)e(\tau)d\tau.
\end{aligned}
\end{equation}
Taking the Laplace transform of \eqref{d.10n}, we have
\begin{equation}\label{d.11}
\begin{aligned}
\hat{\sigma}(s)=s\hat{G}(s)\hat{e}(s).
\end{aligned}
\end{equation}
Comparing \eqref{d.9} with \eqref{d.11}, we have
\begin{equation}\label{d.12}
\begin{aligned}
\hat{G}(s)=\frac{\eta_3+\eta_2\eta_3 c_2^{-1}s}{1+B_1s+B_2s^2}=\frac{\eta_3+\eta_2\eta_3 c_2^{-1}s}{B_2(s+r_1)(s+r_2)},
\end{aligned}
\end{equation}
where $r_1=\frac{B_1+\sqrt{D}}{2B_2}>0$ and $r_2=\frac{B_1-\sqrt{D}}{2B_2}>0$ with $D=B_1^2-4B_2>0$.
We note that the constants of the model satisfy
\begin{equation}\label{d.13}
\begin{aligned}
\eta_2c_2^{-1}r_2\leq 1 \leq \eta_2c_2^{-1}r_1.
\end{aligned}
\end{equation}
Under the condition \eqref{d.13}, we have from \eqref{d.12} that
\begin{equation}\label{d.14}
\begin{aligned}
B_2\hat{G}(s)=\frac{b_1}{s+r_1}+\frac{b_2}{s+r_2},
\end{aligned}
\end{equation}
where $b_1\geq 0$ and $b_2\geq 0$.
Taking the inverse Laplace transform of \eqref{d.14} with respect to $s$, we have
\begin{equation}\label{d.15}
\begin{aligned}
G(t)=b_1e^{-r_1t}+b_2e^{-r_2t}
\end{aligned}
\end{equation}
with
$$
b_1=\frac{\eta_2\eta_3 c_2^{-1}r_1-\eta_3}{B_2(r_1-r_2)}=\frac{\eta_2\eta_3 c_2^{-1}r_1-\eta_3}{\sqrt{D}},\,\,b_2=\frac{\eta_3-\eta_2\eta_3 c_2^{-1}r_2}{B_2(r_1-r_2)}=\frac{\eta_3-\eta_2\eta_3 c_2^{-1}r_2}{\sqrt{D}},
$$
which satisfies all the conditions of Assumption 1.2 except (iv)'. The same can be said about the extended Burgers model which is obtained by connecting Burgers model in series.

The one space dimensional extended Burgers model is obtained by connecting Burgers model in series. If each Burgers model is such a mentioned Burgers model, the kernel of the extended Burgers model can satisfy Assumption \ref{exponential} except (iv)'. If we connect one spring to this model, then \ref{exponential}, (iv)' is satisfied. 

 Finally, we would like to point out our paper \cite{DKLN} which discusses about the property of solutions of initial boundary value problem with mixed type boundary condition for the extended Maxwell system and extended SLS system, and also our paper \cite{EBM} which discusses about the derivation of the relaxation tensor for the higher space-dimensional anisotropic extended Burgers model will help the readers to understand more about these systems. Especially for the higher space dimensional extended Burgers model, one can find the property of the relaxation tensor and the decaying property of solutions of the associated VID system are the same as for the one space dimensional case.

\section{Modification of the proof of Theorem~2.2 in \cite{D1}} \label{App:modification}

The initial boundary value problem (\ref{mixed prob}) has a homogeneous mixed type boundary condition which consists of the homogeneous Dirichlet boundary condition over $\Gamma_D$ and the homogeneous boundary condition over $\Gamma_N$ with the boundary operator $\mathcal{T}+\partial_t$, where $\mathcal{T}$ is the traction operator. In \cite{D1}, the author applied his method to analyze the unique solvability of the initial boundary value problem for VID system with the homogeneous Dirichlet boundary condition over the whole $\partial\Omega$. He commented in the footnote $9$ on page 566 of his paper by saying ``More general classes of homogeneous boundary conditions can be considered alternatively with slight modifications''. Indeed we need a slight modification of the proof for Theorem 2.2, which will be given in this section.

Since allowing to have a source term is important in the proof of Theorem 2.2 in \cite{D1}, we put a source term $f$ to the VID system and consider \eqref{mixed prob} as a Cauchy problem for an abstract Volterra equation (1.2), (1.3) given on page 554 of \cite{D1}.

To begin with, recall  $H_+:=\{\psi\in H^1(\Omega): \psi=0\,\text{on}\,\Gamma_D\}$, $H_0:=L^2(\Omega)$, and the dual space $H_-$ of $H_+$ via the inner product $\ang<\,,\,>$ of $H_0$.

In order to avoid any confusion, we denote our $\rho$, $C=(C_{ijkl}(\cdot))$ and $G(t)=(G_{ijkl}(\cdot,t))$ by $\rho_0$, $C_0$ and $G_0(t)$, respectively. Then, $\rho$, $C$ and $G(t)$ are $\rho=\rho_0 I$, $(Cu)(\cdot,t)=\nabla\cdot(C_0:\nabla u(\cdot,t))$ and $\int_0^t G(t-\tau)u(\tau)d\tau=\int_0^t\nabla\cdot(G_0(t-\tau):\nabla u(\cdot,\tau))d\tau$, respectively. By using the test function $(t-t_0)\dot v$ with a fix $t_0\in(0,T]$ and any $v\in\mathcal{E}_{t_0}:=\{v(t)\in C^\infty([0,t_0]; H_+,\,v(0)=0\}$ to test the VID system, we derive
$$\mathcal{B}(u,v)=\mathcal{D}(f,v)-t_0\ang<\rho\dot u(0), \dot v(0)>,$$
where
$$
\left\{
\begin{array}{ll}
\mathcal{B}(u,v)=\int_0^{t_0}(t-t_0)[\ang<\rho\dot u(t),\ddot v(t)>-\ang<C_0:\nabla u(t),\nabla\dot v(t)>\\
\\
\quad\qquad\qquad+\ang<\int_0^t G_0(t-\tau):\nabla u(\tau)\, d\tau,\nabla\dot v(t)>-s\ang<\dot u(t),\dot v(t)>_{\Gamma_N}]dt\\
\\
\quad\qquad\qquad+\int_0^t \ang<\rho\dot u(t),\dot v(t)>dt,\\
\\
\mathcal{D}(f,v)=-\int_0^{t_0}(t-t_0)[\ang<f^{(1)}(t),\dot v(t)>-\ang<\dot f^{(2)}(t),v(t)>]dt+\int_0^{t_0}\ang<f^{(2)}(t),v(t)>dt,
\end{array}
\right.
$$
where $\ang<\,,\,>_{\Gamma_N}$ is the $L^2(\Gamma_N)$-inner product.
Integrating by parts with respect to $t$, we have
$$
\begin{array}{ll}
\mathcal{B}(v,v)=\frac{1}{2}\int_0^{t_0}\ang<\rho\dot v(t),\dot v(t)>dt+\frac{1}{2}\int_0^{t_0}\ang<C_0;\nabla v(t),\nabla v(t)>dt\\
\\
\qquad-\int_0^{t_0}(t-t_0)[\ang<G_0(0):\nabla v(t),\nabla v(t)>+\ang<\int_0^t\dot G_0(t-\tau):\nabla v(\tau),\nabla v(t)>d\tau]dt\\
\\
\qquad+\int_0^{t_0}\ang<\int_0^t G_0(t-\tau):\nabla v(\tau)d\tau,\nabla v(t)>dt\\
\\
\qquad-s\int_0^{t_0}(t-t_0)\ang<\dot v(t),\dot v(t)>_{\Gamma_N}dt+\frac{1}{2}t_0\ang<\dot v(0),\dot v(0)>.
\end{array}
$$
The final modification is the definition of the inner product $(v,w)_{\Phi}$ for $\mathcal{E}_{t_0}$. We define it as
$$
\begin{array}{ll}
(v,w)_\Phi=&\int_0^{t_0}\{\ang<\dot v(t),\dot w(t)>+<v(t),w(t)>_{H_+}\}dt-s\int_0^{t_0}(t-t_0)\ang<\dot v(t),\dot w(t)>_{\Gamma_N}dt\\
&+t_0\ang<\dot v(0),\dot w(0)>.
\end{array}
$$

\end{document}